\documentclass[11pt, reqno]{amsart}

\usepackage[usenames,dvipsnames]{color}

\usepackage[OT2, T1]{fontenc}
\usepackage{url}
\usepackage{amsmath}
\usepackage{array}
\usepackage{graphicx}
\usepackage{amssymb}
\usepackage{amsthm}
\definecolor{link}{RGB}{11,0,128}
\usepackage[colorlinks=true,citecolor=NavyBlue,linkcolor=OliveGreen,urlcolor=link]{hyperref}
\usepackage{hyperref}
\usepackage{paralist}
\usepackage{colonequals}			
\usepackage{color}
\usepackage{mathabx}			
\usepackage{amscd}
\usepackage[all,cmtip]{xy}			
\usepackage{sseq}				
\usepackage{verbatim}
\usepackage{parskip}
\usepackage{microtype}
\usepackage[in]{fullpage}
\usepackage{mathrsfs} 			
\usepackage{cleveref}
\usepackage{enumitem}
\usepackage{moreenum}			
\usepackage[alphabetic,lite]{amsrefs} 	
\usepackage{stmaryrd}			
\usepackage[foot]{amsaddr}
\usepackage{subfiles}
\usepackage{ragged2e}
\usepackage{stackengine}

\newcommand{\gA}{\alpha}

\newcommand{\bA}{\mathbb{A}}

\newcommand{\bF}{\mathbb{F}}
\newcommand{\bG}{\mathbb{G}}

\newcommand{\bP}{\mathbb{P}}

\newcommand{\bZ}{\mathbb{Z}}

\newcommand{\bbB}{\mathbf{B}}

\newcommand{\cO}{\mathcal{O}}

\newcommand{\cX}{\mathcal{X}}

\newcommand{\fm}{\mathfrak{m}}
\newcommand{\fn}{\mathfrak{n}}

\newcommand{\fp}{\mathfrak{p}}

\newcommand{\sE}{\mathscr{E}}
\newcommand{\sF}{\mathscr{F}}
\newcommand{\sG}{\mathscr{G}}

\newcommand{\sM}{\mathscr{M}}

\newcommand{\sO}{\mathscr{O}}

\newcommand{\sR}{\mathscr{R}}

\newcommand{\sV}{\mathscr{V}}

\newcommand{\sY}{\mathscr{Y}}
\newcommand{\sZ}{\mathscr{Z}}

\DeclareMathOperator{\Aut}{Aut}		
\DeclareMathOperator{\Bl}{Bl}			
\DeclareMathOperator{\Bun}{Bun}			
\DeclareMathOperator{\Corad}{Corad}	
\DeclareMathOperator{\Ext}{Ext}		
\DeclareMathOperator{\Frac}{Frac}		
\DeclareMathOperator{\Hom}{Hom}		
\DeclareMathOperator{\im}{Im}			
\DeclareMathOperator{\Isom}{Isom}		
\DeclareMathOperator{\Ker}{Ker}		
\DeclareMathOperator{\rad}{rad}		
\DeclareMathOperator{\Rad}{Rad}		
\DeclareMathOperator{\Res}{Res}		
\DeclareSymbolFont{cyrletters}{OT2}{wncyr}{m}{n}
\DeclareMathSymbol{\Sha}{\mathalpha}{cyrletters}{"58}	
\DeclareMathOperator{\Spec}{Spec}		

\newcommand{\ad}{\mathrm{ad}}			

\newcommand{\ce}{\colonequals}

\newcommand{\der}{\mathrm{der}}		

\newcommand{\eps}{\varepsilon}
\newcommand{\fppf}{\mathrm{fppf}}		

\newcommand{\gp}{{\mathrm{gp}}}		
\newcommand{\hra}{\hookrightarrow}
\renewcommand{\i}{^{-1}}

\newcommand{\isomto}{\overset{\sim}{\longrightarrow}}

\newcommand{\llb}{\llbracket}			
\newcommand{\llp}{(\!(}			
	
\newcommand{\ov}{\overline}

\newcommand{\ra}{\rightarrow}

\newcommand{\rrb}{\rrbracket}			
\newcommand{\rrp}{)\!)}			

\newcommand{\surjects}{\twoheadrightarrow}
\newcommand{\tensor}{\otimes} 			
\newcommand{\tors}{\mathrm{tors}}		

\newcommand{\wt}{\widetilde}
\newcommand{\xra}{\xrightarrow}

\providecommand{\up}[1]{{\upshape(}#1{\upshape)}}
\providecommand{\uref}[1]{{\upshape\ref{#1}}}
\providecommand{\uS}{{\upshape\S}}

\renewcommand{\b}{\textbf}
\providecommand{\ucolon}{{\upshape:} }
\providecommand{\uscolon}{{\upshape;} }

\newcommand{\brems}{\begin{rems} \hfill \begin{enumerate}[label=\b{\thenumberingbase.},ref=\thenumberingbase]}
\newcommand{\remi}{\addtocounter{numberingbase}{1} \item}
\newcommand{\erems}{\end{enumerate} \end{rems}}
\newcommand{\begs}{\begin{egs} \hfill \begin{enumerate}[label=\b{\thenumberingbase.},ref=\thenumberingbase]}

\newcommand{\eegs}{\end{enumerate} \end{egs}}
\newcommand{\m}{\item}
\newcommand{\bsm}{\begin{smallmatrix}}
\newcommand{\esm}{\end{smallmatrix}}
\newcommand{\blem}{\begin{lemma}}
\newcommand{\elem}{\end{lemma}}

\newcommand{\bconj}{\begin{conj}}
\newcommand{\econj}{\end{conj}}
\newcommand{\bprob}{\begin{Problem}}
\newcommand{\eprob}{\end{Problem}}

\newcommand{\bq}{\begin{Q}}
\newcommand{\eq}{\end{Q}}
\newcommand{\benum}{\begin{enumerate}[label={{\upshape(\alph*)}}]}
\newcommand{\benuma}{\begin{enumerate}[label={{\upshape(\arabic*)}}]}
\newcommand{\benumb}{\begin{enumerate}[label={{\upshape\b{\arabic*.}}}]}
\newcommand{\benumr}{\begin{enumerate}[label={{\upshape(\roman*)}}]}
\newcommand{\eenum}{\end{enumerate}}
\newcommand{\bitem}{\begin{itemize}}
\newcommand{\eitem}{\end{itemize}}
\newcommand{\bc}{}
\newcommand{\bd}{\begin{defn}}
\newcommand{\ed}{\end{defn}}

\newcommand{\beg}{\begin{eg}}
\newcommand{\eeg}{\end{eg}}

\newcommand{\bcl}{\begin{claim}}
\newcommand{\ecl}{\end{claim}}

\newcommand{\x}{\text}

\newcommand{\q}{\quad}
\providecommand{\qxq}[1]{\quad\text{#1}\quad}

\newcommand{\qq}{\quad\quad}

\newcommand{\tst}{\textstyle}
\newcommand{\ba}{\begin{aligned}}
\newcommand{\ea}{\end{aligned}}
\newcommand{\be}{\begin{equation}}
\newcommand{\ee}{\end{equation}}
\newcommand{\bpf}{\begin{proof}}
\newcommand{\epf}{\end{proof}}
\newcommand{\bthm}{\begin{thm}}
\newcommand{\ethm}{\end{thm}}
\newcommand{\bprop}{\begin{prop}}
\newcommand{\eprop}{\end{prop}}
\newcommand{\bcor}{\begin{cor}}
\newcommand{\ecor}{\end{cor}}
\newcommand{\brem}{\begin{rem}}
\newcommand{\erem}{\end{rem}}

\usepackage{aliascnt}
\newaliascnt{numberingbase}{subsection}
\numberwithin{equation}{section}

\newtheoremstyle{thms}{0.5em}{0.5em}{\itshape}{}{\bfseries}{.}{ }{}
\theoremstyle{thms}
\newtheorem{conj}[numberingbase]{Conjecture}

\newtheorem{cor}[numberingbase]{Corollary}

\newtheorem{lemma}[numberingbase]{Lemma}

\newtheorem{prop}[numberingbase]{Proposition}

\newtheorem{Q}[numberingbase]{Question}
\newtheorem{thm}[numberingbase]{Theorem}

\newtheoremstyle{claims}{0.5em}{0.5em}{}{}{\itshape}{.}{ }{}
\theoremstyle{claims}
\newtheorem{claim}[numberingbase]{Claim}

\newtheoremstyle{defs}{0.5em}{0.5em}{}{}{\bfseries}{.}{ }{}
\theoremstyle{defs}
\newtheorem{defn}[numberingbase]{Definition}

\newtheorem{eg}[numberingbase]{Example}
\newtheorem*{egs}{Examples}
\newtheorem{rem}[numberingbase]{Remark}

\newtheorem*{rems}{Remarks}

\Crefname{claim}{Claim}{Claims}
\Crefname{bclaim}{Claim}{Claims}
\Crefname{sublemma}{Lemma}{Lemmas}
\Crefname{conj}{Conjecture}{Conjectures}
\Crefname{cor}{Corollary}{Corollaries}
\Crefname{defn}{Definition}{Definitions}
\Crefname{eg}{Example}{Examples}
\Crefname{prop}{Proposition}{Propositions} 
\Crefname{Q}{Question}{Questions}
\Crefname{rem}{Remark}{Remarks}
\Crefname{thm}{Theorem}{Theorems}
\Crefname{Theorem}{Theorem}{Theorems}
\Crefname{variant}{Variant}{Variants}
\Crefname{caution}{Caution}{Cautions}

\theoremstyle{thms}
\newtheorem{thm-tweak}[subsection]{Theorem}
\Crefname{thm-tweak}{Theorem}{Theorems}
\newtheorem{lemma-tweak}[subsection]{Lemma}
\Crefname{lemma-tweak}{Lemma}{Lemmas}
\newtheorem{cor-tweak}[subsection]{Corollary}
\Crefname{cor-tweak}{Corollary}{Corollaries}
\newtheorem{prop-tweak}[subsection]{Proposition}
\Crefname{prop-tweak}{Proposition}{Propositions} 
\newtheorem{conj-tweak}[subsection]{Conjecture}
\Crefname{conj-tweak}{Conjecture}{Conjectures} 
\newtheorem{q-tweak}[subsection]{Question}
\Crefname{q-tweak}{Question}{Questions} 

\theoremstyle{defs}
\newtheorem{defn-tweak}[subsection]{Definition}
\Crefname{defn-tweak}{Definition}{Definitions}
\newtheorem{eg-tweak}[subsection]{Example}
\Crefname{eg-tweak}{Example}{Examples}
\newtheorem*{rems-tweak}{Remarks}
\newtheorem{rem-tweak}[subsection]{Remark}
\Crefname{rem-tweak}{Remark}{Remarks}

\newtheoremstyle{subsection-tweak}
   {2pt}
   {3pt}%
   {}
   {}%
   {\bfseries}
   {}%
   {.5em}
   {\thmnumber{\@{#1}{}\@{#2}.}%
    \thmnote{~{\bfseries#3.}}}    
    
\theoremstyle{subsection-tweak}
\newtheorem{pp}[numberingbase]{}
\newcommand{\bpp}{\begin{pp}}
\newcommand{\epp}{\end{pp}}

\theoremstyle{subsection-tweak}
\newtheorem{pp-tweak}[subsection]{}

\makeatletter
\def\@tocline#1#2#3#4#5#6#7{
    \begingroup 
    \@ifempty{#4}{}{}

    \parindent\z@ \leftskip#3\relax \advance\leftskip\@tempdima\relax
    #5\hskip-\@tempdima
      \ifcase #1
       \or\or \hskip 2em \or \hskip 1em \else \hskip 3em \fi%
      #6\nobreak\relax
    \dotfill\hbox to\@pnumwidth{\@tocpagenum{#7}}\par
    \nobreak
    \endgroup
 }
 \def\l@section{\@tocline{1}{0pt}{1pc}{}{}}

\renewcommand{\tocsection}[3]{%
  \indentlabel{\@ifnotempty{#2}{\makebox[1.3em][l]{%
    \ignorespaces#1 \bfseries{#2}.\hfill}}}\bfseries{#3}
    \vspace{-5pt}}

\renewcommand{\tocsubsection}[3]{%
  \indentlabel{\@ifnotempty{#2}{\hspace*{-0.5em}\makebox[2.1em][l]{%
    \ignorespaces#1#2.\hfill}}}#3
    \vspace{-5pt}}
   
\makeatother 

\makeatletter 
\newcommand\appendix@section[1]{%
  \refstepcounter{section}%
  \orig@section*{Appendix \@Alph\c@section. #1}%
}
\let\orig@section\section
\g@addto@macro\appendix{\let\section\appendix@section}
\makeatother

\makeatletter
\@namedef{subjclassname@2020}{%
  \textup{2020} Mathematics Subject Classification}
\makeatother

\author{K\k{e}stutis \v{C}esnavi\v{c}ius}
\author{Roman Fedorov}
\address{Sorbonne Universit\'e, Universit\'e Paris Cit\'e, CNRS, IMJ-PRG, F-75005 Paris, France}
\address{University of Pittsburgh, 301 Thackeray Hall, Pittsburgh, PA 15260, USA}
\address{Max-Planck-Institut f\"ur Mathematik, Vivatsgasse 7, 53111 Bonn, Germany}
\email{kestutis@imj-prg.fr, fedorov@pitt.edu}

\date{\today}

\begin{document}


\title{Unramified Grothendieck--Serre for isotropic groups}

\maketitle

\begin{abstract} 
The Grothendieck--Serre conjecture predicts that every generically trivial torsor under a reductive group $G$ over a regular semilocal ring $R$ is trivial. We establish this for unramified $R$ granted that $G^{\mathrm{ad}}$ is totally isotropic, that is, has a ``maximally transversal'' parabolic $R$-subgroup. We also use purity for the Brauer group to reduce the conjecture for unramified $R$ to simply connected $G$---a much less direct such reduction of Panin had been a step in solving the equal characteristic case of Grothendieck--Serre. We base the group-theoretic aspects of our arguments on the geometry of the stack $\Bun_G$, instead of the affine Grassmannian used previously, and we quickly reprove the crucial weak $\bP^1$-invariance input: for any reductive group $H$ over a semilocal ring $A$, every $H$-torsor $\sE$ on $\bP^1_A$ satisfies $\sE|_{\{t = 0\}} \simeq \sE|_{\{t = \infty\}}$. For the geometric aspects, we develop reembedding and excision techniques for relative curves with finiteness weakened to quasi-finiteness, thus overcoming a known obstacle in mixed characteristic, and show that every generically trivial torsor over $R$ under a totally isotropic $G$ trivializes over every affine open of $\Spec(R) \setminus Z$ for some closed $Z$ of codimension $\ge 2$.

 \end{abstract}

\hypersetup{
    linktoc=page,     
}
\tableofcontents

\section{The unramified totally isotropic case of the Grothendieck--Serre conjecture}

In this article, we solve a case of the following conjecture of Grothendieck and Serre \cite{Ser58b}*{page 31, remarque}, \cite{Gro58}*{pages~26--27, remarques 3}, \cite{Gro68b}*{remarques 1.11 a)}  about triviality of torsors.

\bconj[Grothendieck--Serre] \label{conj:GS}
For a reductive group scheme over a regular semilocal ring $R$, no nontrivial $G$-torsor over $R$ trivializes over the total ring of fractions $K \ce \Frac(R)$, that is,
\[
\Ker(H^1(R, G) \ra H^1(K, G)) = \{ *\}. 
\]
\econj


Torsors occur naturally in many contexts, for instance, in studying conjugacy of sections. For conjugacy problems, \Cref{conj:GS} predicts that conjugacy over $K$ of sections over $R$ implies conjugacy over $R$, granted that the centralizer group schemes are reductive and fiberwise connected.

The Grothendieck--Serre conjecture is a nonabelian avatar of Gersten injectivity conjectures for various abelian cohomology theories of motivic flavor. Indeed, one may hope that $H^1(R, G)$ could be described in terms of abelian cohomological invariants in the style of \cite{Ser94}*{sections 6--10}, at which point \Cref{conj:GS} would follow from these abelian counterparts. Such an approach is firmly out of reach of available technology, but it is plausible that it could eventually be reversed, namely, that \Cref{conj:GS} may eventually be used to describe $H^1(R, G)$ by abelian cohomological invariants.

We settle the Grothendieck--Serre conjecture in the case when the regular ring $R$ is unramified and the group $G$ is such that its adjoint quotient $G^\ad$ has no anisotropic factors.

\bthm[\Cref{thm:tot-iso}] \label{thm:main}
Let $R$ be a Noetherian semilocal  ring that is flat and geometrically~regular\footnote{We recall from \cite{SP}*{Definition \href{https://stacks.math.columbia.edu/tag/0382}{0382}} that the geometric regularity assumption means that $R \tensor_k k'$ is a regular ring for every finite extension $k'$ of some residue field $k$ of $\cO$. By Popescu theorem \cite{SP}*{Theorem \href{https://stacks.math.columbia.edu/tag/07GC}{07GC}}, it is equivalent to require that our regular semilocal $R$ be a filtered direct limit of smooth $\cO$-algebras.} over a Dedekind ring $\cO$, let $K \ce \Frac(R)$ be its fraction ring, and let $G$ be a reductive~$R$-group such that $G^{\mathrm{ad}}$ is totally isotropic \up{see \eqref{eqn:tot-isot}}. No nontrivial $G$-torsor over $R$ trivializes over $K$, that is, 
\[
\q \Ker(H^1(R, G) \ra H^1(K, G)) = \{ *\}.
\]
\ethm


The following are the cases in which the Grothendieck--Serre conjecture has been established.
\benumr
\m \label{m:known-i}
In equal characteristic, that is, when $\cO$ in \Cref{thm:main} is a field, the Grothendieck--Serre conjecture was settled by Fedorov--Panin \cite{FP15} and Panin \cite{Pan20a}, with simplifications in \cite{Fed22a} and significant special cases obtained in prior works \cite{Oja80}, \cite{CTO92}, \cite{Rag94}, \cite{PS97}, \cite{Zai00}, \cite{OP01}, \cite{OPZ04}, \cite{Pan05}, \cite{Zai05}, \cite{PPS09}, \cite{PS09}, \cite{Che10}, \cite{PSV15},~\cite{Pan20b}; see also \cite{Pan22a} for a variant beyond connected reductive groups. 

\m \label{m:known-ii}
For regular semilocal $R$ that are unramified, more precisely, that are as in \Cref{thm:main}, the Grothendieck--Serre conjecture has been established for quasi-split $G$ in \cite{split-unramified} (with a prior more restrictive case in \cite{Fed21a}) and for $G$ that descend to reductive $\cO$-groups in \cite{GL24a} (with subcases of this constant case already in \cite{Pan19c}, \cite{GP23}). For further variants with, more generally, $\cO$ a semilocal Pr\"ufer ring of dimension $\le 1$, see \cite{GL24a}, \cite{GL24}*{Theorem~8.1}, and \cite{Kun23}*{Theorem A on page 24} (the latter with $\cO$ a valuation ring of dimension $\le 1$).  

\m \label{m:known-iii}
The conjecture is known in the case when $R$ is of dimension $\le 1$ by \cite{Guo22a} that built on prior \cite{Nis82} and \cite{Nis84} (with special cases in \cite{Har67}, \cite{BB70}, \cite{BT87}, \cite{PS16}, \cite{BvG14}, \cite{BFF17}, \cite{BFFH19}, and valuation ring variants in \cite{Guo22b} and \cite{GL24a}*{Appendix A}). This one-dimensional case implies the case when $R$ is Henselian, see \cite{CTS79}*{assertion~6.6.1}.

\m \label{m:known-iv}
The case when $G$ is a torus was settled by Colliot-Th\'{e}l\`{e}ne and Sansuc in \cite{CTS78}, \cite{CTS87}. 

\m \label{m:known-v}
Sporadic cases with either $G$ or $R$ of specific form were settled in \cite{Gro68b}*{remarques 1.11 a)}, \cite{Oja82}, \cite{Nis89}, \cite{BFFP22}, \cite{Fir23}, \cite{Pan22b}.
\eenum
 
For arguing \Cref{thm:main}, we only use the $1$-dimensional case \ref{m:known-iii}, but not any of the other cases.

Throughout the works above, there are broadly two approaches to the Grothendieck--Serre conjecture: 
\bitem
\m
the geometric approach, which was pioneered by Colliot-Th\'el\`ene--Ojanguren \cite{CTO92} and then developed much further in the works that culminated in the results \ref{m:known-i}--\ref{m:known-ii}; and


\m
the group-theoretic approach, prevalent in \ref{m:known-iii}--\ref{m:known-v} and based on analyzing the structure of~$G$.
\eitem
The group-theoretic approach appeared earlier, and its ideas and results later fed into the  geometric approach, which analyzes the interaction of the geometry of $R$ with the properties of $G$. Given a generically trivial $G$-torsor $E$ over $R$, the gist of the geometric approach is to explicate the geometry of $R$ via presentation lemmas of Gabber--Quillen type and to combine them with patching arguments to eventually produce a $G$-torsor $\sE$ over $\bP^1_R$ such that $\sE|_{\{ t = 0\}} \simeq E$ and $\sE|_{\{t = \infty\}}$ is trivial. 
On the other hand, results rooted in the geometry of the algebraic stack $\Bun_G$ parametrizing $G$-torsors over the relative projective line imply that every family of $G$-torsors over $\bP^1_R$ is $R$-sectionwise constant, in particular, that $\sE|_{\{t = 0\}} \simeq \sE|_{\{t = \infty\}}$, see \Cref{thm:section-const} below or \cite{PS25}*{Theorem 1.2}. Taken together, this means that $E$ is trivial.

In this article, we develop the geometric approach further, the following being our main novelties.
\benuma
\m \label{m:1}
In comparison to equal characteristic, the main complication in the unramified mixed characteristic case of the Grothendieck--Serre conjecture is that the base $\cO$ of the projection that we have no flexibility to ``move'' is now one-dimensional, which makes us lose one dimension in geometric arguments. For instance, to start the geometric approach we now have to build a closed $Z \subset \Spec R$ \emph{of codimension $\ge 2$} away from which our generically trivial $G$-torsor $E$ over $R$ is ``simpler,'' whereas in equal characteristic (when $\cO$ was a field) codimension $\ge 1$ sufficed and was straight-forward to arrange from generic triviality. In \S\ref{sec:codim-2}, we bypass this problem: for any $G$ and $E$, in \Cref{prop:build-Z}, we build an open $V \subset \bP^1_R$ containing both $\bP^1_{\Spec(R) \setminus Z}$ for some closed $Z \subset \Spec(R)$ of codimension $\ge 2$ and the sections $\{t = 0\}$ and $\{t = \infty\}$, as well as a $G$-torsor $\sE$ over $V$ such that $\sE|_{\{t = 0\}} \simeq E$ and $\sE|_{\{t = \infty\}}$ is trivial. 

Consequently, $E$ becomes ``simpler'' over $\Spec(R) \setminus Z$ in the sense that it fits into a family of $G$-torsors over $\bP^1_{\Spec(R) \setminus Z}$ with a trivial fiber at infinity. For $G$ with $G^{\mathrm{ad}}$ totally isotropic, this already implies that $E$ trivializes over every affine $(\Spec(R) \setminus Z)$-scheme, see \Cref{thm:Horrocks}. 

To build $V$, we use a quasi-finite version of the presentation lemma and find a way to carry out the subsequent reembedding techniques with finiteness weakened to quasi-finiteness. In contrast, 
 building the desired $Z$ of codimension $\ge 2$ was simpler in \cite{split-unramified}: it sufficed to combine the quasi-splitness assumption made there with the valuative criterion of properness. 

\m \label{m:2}
We take advantage of our $\sE$ over $V$ as in \ref{m:1} in several different (and disjoint) ways.

Firstly, in \S\ref{sec:tot-iso}, we use our $\sE$ and $V$ to carry out the geometric approach in full for totally isotropic $G$: we settle the unramified case of the Grothendieck--Serre conjecture for such $G$ in \Cref{thm:tot-iso}. Roughly, $\sE$ and $V$ serve as witnesses of $E$ being ``simpler'' over $\Spec(R) \setminus Z$, and we carry them along the steps of the geometric approach to eventually build a $G$-torsor $\sF$ over $\bP^1_R$ (unrelated to $\sE$) such that $\sF|_{\{t = 0\}} \simeq E$ and $\sF|_{\{t = \infty\}}$ is trivial. The $R$-sectionwise constancy of families of $G$-torsors over $\bP^1_R$ applied to $\sF$ then implies the triviality of $E$. 

A crucial novel aspect of our implementation of the geometric approach is to carry along not only $Z$, but also a closed $Y \subset \Spec(R)$ of codimension $1$ containing it such that $E|_{\Spec(R) \setminus Y}$ is trivial: $Y$ is important for mitigating the loss of applicability of the excision lemma for unipotent torsors \cite{split-unramified}*{Lemma 7.2 (b)} to pass to $\bP^1_R$ in our setting. Relatedly, in \Cref{lem:pres-lem} we generalize the mixed characteristic presentation lemma to track both $Y$ and $Z$.

Secondly, in \S\ref{sec:reduce-ss}, we combine the existence of $\sE$ with the purity for the Brauer group (see \cite{brauer-purity}) and constancy for multiplicative group gerbes over $\bP^1_R$ (see \Cref{lem:gerbeM}) to quickly reduce the unramified case of  Grothendieck--Serre to simply connected groups. This method is new even in equicharacteristic, where the corresponding result was the main goal of \cite{Pan20b}.

\m \label{m:3}
For studying $G$-torsors over a relative $\bP^1$, we base our arguments on the geometry of the algebraic moduli stack $\Bun_G$ parametrizing such torsors. This replaces affine Grassmannian inputs used in previous works starting with \cite{FP15} and leads to clean, simple, broadly useful geometric statements about $\Bun_G$ recorded in \S\ref{sec:BunG}, for instance, \Cref{prop:const} or \Cref{thm:section-const}.

\eenum

Even though we limit ourselves to the totally isotropic unramified case, our results may also reach most types of anisotropic reductive $G$ over an unramified regular semilocal $R$ as follows. First of all, by passing to the simply connected case via \Cref{prop:pass-to-sc} and decomposing into factors, we may harmlessly assume that $G$ has simple fibers. The main idea then comes from observing that if $G \hra \wt{G}$ is an inclusion of a factor of a  Levi subgroup of a larger reductive $R$-group $\wt{G}$, then 
\[
H^1(R, G) \hra H^1(R, \wt{G}) \qxq{and} H^1(K, G) \hra H^1(K, \wt{G}),
\]
see, for instance \cite{torsors-regular}*{equation (1.3.5.2)}. This reduces the Grothendieck--Serre conjecture for $G$ to that for $\wt{G}$; however, the latter is isotropic, so \Cref{thm:main} applies to it. The focus then shifts to realizing $G$ inside some $\wt{G}$ in this way. Overall, this type of approach to anisotropic groups was explored in \cite{PPS09} in equal characteristic, but one may amplify it further by first combining techniques of \S\ref{sec:codim-2} with ideas from \cite{Pan20b} to obtain the flexibility of varying $G$ in isogenies or even passing to studying generically isomorphic adjoint $R$-groups instead of torsors. Nevertheless, even though we could reach most types of anisotropic $G$ in this way, types such as $F_4$ or $E_8$ never occur as Levis of larger reductive groups and seem too large to treat directly, which signals the need of other ideas for arguing the remaining anisotropic case for unramified $R$ in a clean conceptual way.




\bpp[Notation and conventions] \label{pp:conv}
For a field $k$, we let $\ov{k}$ denote its algebraic closure. For a point $s$ of a scheme $S$ (resp.,~a prime ideal $\fp$ of a ring $R$), we let $k_s$ (resp.,~$k_\fp$) denote its residue field viewed as an algebra over $S$ (resp.,~over $R$). We let $\Frac(-)$ denote both the total ring of fractions of a ring and the function field of an integral scheme, depending on the context. 

When it comes to reductive groups, we follow SGA 3, in particular, a reductive group over a scheme $S$ is a smooth, affine $S$-group scheme whose geometric fibers are connected reductive groups, see \cite{SGA3IIInew}*{expos\'{e} XIX, d\'efinition 2.7}. See also \cite{torsors-regular}*{Section 1.3} for a review of basic reductive group notions and notations that we use freely. In particular, we write $G^\der$ for the derived subgroup of a reductive group scheme $G$ and we write $H^{\mathrm{sc}}$ for the simply connected cover of a semisimple group scheme $H$ (see \emph{loc.~cit.}~for a review). Similarly to \cite{split-unramified}*{Definition~8.1} (or \cite{torsors-regular}*{Section~1.3.6}), a semisimple $S$-group $G$ is \emph{totally isotropic} if in the canonical decomposition
\be \label{eqn:tot-isot}
\tst G^\ad \cong \prod_{i \in \{A_n, B_n, \dotsc, G_2\}} \Res_{S_i/S}(G_i)
\ee
of \cite{SGA3IIInew}*{expos\'{e} XXIV, proposition 5.10 (i)}, in which $i$ ranges over the types of connected Dynkin diagrams, $S_i$ is a finite \'{e}tale $S$-scheme, and $G_i$ is an adjoint semisimple $S_i$-group with simple geometric fibers of type $i$, Zariski locally on $S$ each $G_i$ has a parabolic $S_i$-subgroup that contains no $S_i$-fiber of $G_i$; intuitively, this amounts to requiring that Zariski locally on $S$ the group $G$ itself contain a proper (relative to each factor) parabolic subgroup. 

We say that a reductive $S$-group $G$ is \emph{simple} if it is semisimple and the Dynkin diagrams of its geometric $S$-fibers are all connected (some authors call such groups absolutely almost simple because even in the case when $S$ is a geometric point, $G$ may still have nontrivial finite central subgroups). 
\epp

\subsection*{Acknowledgements} 
We thank the referees for helpful comments and suggestions. We thank Alexis Bouthier, Ofer Gabber, Philippe Gille, Ning Guo, Arnab Kundu, Ivan Panin, and David Rydh for helpful conversations and correspondence. We thank the Max Planck Institute for Mathematics and the Institute for Advanced Study for excellent working conditions; significant portions of this project were carried out while enjoying the hospitality of these institutions. This project has received funding from the European Research Council (ERC) under the European Union's Horizon 2020 research and innovation programme (grant agreement No.~851146). This project is based upon work supported by the National Science Foundation under Grant No.~DMS-1926686 and Grant No.~DMS-2001516.


\section{Lifting to a family of torsors over $\bP^1_R$ away from a closed of codimension $\ge 2$} \label{sec:codim-2}

Our first goal is \Cref{prop:build-Z} below that builds a closed $Z \subset \Spec R$ \emph{of codimension $\ge 2$}, away from which our generically trivial $G$-torsor $E$ over $R$ simplifies. The construction of this $Z$ ultimately hinges on the quasi-finite version of the Gabber--Quillen presentation lemma in mixed characteristic that we establish in \Cref{lem:pres-lem}. This proposition involves a scheme $X$ that is smooth over a semilocal Dedekind domain $\cO$ and a finite set of points of $X$. We start with two lemmas that will allow us to assume that all of these points specialize to closed points of some closed $\cO$-fiber of $X$.

\blem \label{lem:semilocal-limit}
Let $K$ be a field.
\benumr
\m \label{m:SL-i}
The map sending a $1$-dimensional semilocal Dedekind domain $\cO$ with fraction field $K$ to the set of its localizations $\{\cO_\fm\}$ at the maximal ideals $\fm \subset \cO$ gives a bijection between the set of such subrings $\cO$ of $K$ and the set of finite sets of discrete valuations on $K$ \up{encoded by their corresponding valuation rings}\uscolon the inverse of this bijection sends $\{\cO_\fm\}$ to $\bigcap \cO_\fm$.

\m \label{m:SL-ii}
For a subfield $K' \subset K$ and $\cO$ as in \ref{m:SL-i}, the ring $\cO \cap K'$ is a semilocal Dedekind domain with fraction field $K'$.

\m \label{m:SL-iii}
Any $\cO$ as in \ref{m:SL-i} is a filtered direct union of semilocal Dedekind
subdomains whose fraction fields are finitely generated over the prime
field of $K$.
\eenum
\elem

\bpf
Part \ref{m:SL-i} is essentially a restatement of \cite{Mat89}*{Theorem 12.2}. For part \ref{m:SL-ii} it is enough to intersect the equality $\cO=\bigcap_\fm \cO_\fm$
with $K'$ (note that it might happen that $\cO\cap K'=K'$, which is still
a Dedekind ring). Finally, part \ref{m:SL-iii} follows from part \ref{m:SL-ii}.
\epf

\blem \label{lem:Ded-spread}
Let $\cO$ be a semilocal Dedekind domain  whose fraction field $K$ is finitely generated over its prime field, let $X$ be a smooth affine scheme of pure relative dimension $d > 0$ over $\cO$, and let $x_1, \dotsc, x_n \in X$ be finitely many points. There are a semilocal Dedekind domain $\wt{\cO}$ with fraction field $K$ such that $\cO$ is a localization of $\wt{\cO}$ \up{equivalently, $\wt{\cO} \subset \cO$, see Lemma \uref{lem:semilocal-limit}~\uref{m:SL-i}} and a smooth affine $\wt{\cO}$-scheme $\wt{X}$ of pure relative dimension $d > 0$ extending $X$ such that each $x_i$ specializes to a closed point of some closed $\wt{\cO}$-fiber of $\wt{X}$. 
\elem

\bpf
If $\cO$ is a field, then it suffices to take $\wt{\cO} \ce \cO$ and note that any finite type scheme over a field, such as each closure $\overline{\{x_i\}} \subset X$, has a closed point, see \cite{SP}*{Lemma \href{https://stacks.math.columbia.edu/tag/02J6}{02J6}}. In the remaining case when $\cO$ is $1$-dimensional, the same argument shows that each $x_i$ specializes to a closed point of its $\cO$-fiber of $X$, and thus also to the closed point of our sought $\wt{X}$ that we will build later. Thus, we may assume that all the $x_i$ belong to the generic fiber $X_K$, in fact, that they are closed points of $X_K$ (we do not need to worry about the points lying in the closed fibers of $X$ because they will automatically lie in the closed fibers of $\wt{X}$).  

To conclude, we will now follow the ``glue in a new discrete valuation ring'' argument given in the proof of \cite{split-unramified}*{Variant 3.7}. Namely,  since the prime subfield of $K$ is perfect, $K$ is a separable extension of its prime subfield, so \cite{EGAIV4}*{corollaire 17.15.9} and spreading out ensure that $K$ is the fraction field of a domain $A$ that is a smooth algebra either over $\bZ$ or over some $\mathbb{F}_p$. By spreading out $X_K$ and the $x_i$ and replacing $A$ with its localization (that is still a smooth algebra either over $\bZ$ or $\bF_p$), we may assume that $X_K$ extends to a smooth affine $A$-scheme $\cX$ of pure relative dimension $d > 0$ and that each $x_i$ spreads out to an $A$-finite closed subscheme $Z_i \subset \cX$. Since $\cO$ is $1$-dimensional, $K$ is infinite, so $A$ is positive dimensional (if $A$ was zero dimensional, so a field, then, since it is of finite type over $\bZ$, it would have to be a finite field). In particular, since $A$ is  of finite type over $\bZ$, it has infinitely many prime ideals $\fp$ of height $1$, so for some such $\fp$ the discrete valuation subring $A_\fp \subset K$ is different from each $\cO_\fm \subset K$. \Cref{lem:semilocal-limit}~\ref{m:SL-i} ensures that $\cO$ is a localization of the semilocal Dedekind domain $\wt{\cO} \ce A_\fp \cap \cO$ with fraction field $K$ and that, by glueing of $\cX_{A_\fp}$ with $X$, we obtain a smooth affine $\wt{\cO}$-scheme $\wt{X}$ with $\wt{X}_{\wt{\cO}_\fp} \cong \cX_{A_\fp}$ and $\cX_\cO \cong X$. It remains to note that, by construction of the $Z_i$, each $x_i$ specializes to a closed point of the $\fp$-fiber~of~$\wt{X}$. 
\epf

\bprop \label{lem:pres-lem}
Let $X$ be a smooth affine scheme of pure relative dimension $d > 0$ over a Dedekind ring $\cO$, let $x_1, \dotsc, x_n \in X$, and let $Z \subset Y \subset X$ be closed not containing any irreducible component of any $\cO$-fiber of $X$. If either  $Z$ is of codimension $\ge 2$ in $X$ or if $\cO$ is $0$-dimensional, then there are an affine open $X' \subset X$ containing all the $x_i$, an affine open $S \subset \bA^{d - 1}_\cO$, and a smooth $\cO$-morphism $\pi \colon X' \ra S$ such that $Y \cap X'$ is $S$-quasi-finite and $Z \cap X'$ is $S$-finite. 
\eprop
 
 \bpf
 With $Y = Z$, the claim was settled in \cite{split-unramified}*{Proposition~4.1, Remark 4.3}, in fact, it was one of the main technical results of \emph{op.~cit.} We will obtain the general case by similar arguments. 
 
 By localizing at the images of the $x_1, \dotsc, x_n$ in $\Spec \cO$ and then spreading out, we may assume that $\cO$ is semilocal and then, by passing to components, that $\cO$ is a domain. Moreover, by \Cref{lem:semilocal-limit}~\ref{m:SL-iii} and a limit argument, we may assume that the fraction field $K$ of $\cO$ is finitely generated over its prime field. We then enlarge $\cO$ and $X$ as in \Cref{lem:Ded-spread} (and replace $Z$ and $Y$ by their corresponding closures in this larger $X$) to arrange that each $x_i$ that lies in the generic $\cO$-fiber of $X$ has a specialization that lies in some closed $\cO$-fiber of $X$. By replacing such $x_i$ by these specializations, we are therefore left with the case when each $x_i$ lies in some closed $\cO$-fiber of $X$ and is a closed point~of~$X$.
 
At this point, we embed $X$ into an affine space over $\cO$ and form the closure in the corresponding projective space to build an open immersion $X \hra \ov{X}$ into a projective $\cO$-scheme $\ov{X}$, which is flat by \cite{SP}*{Lemma \href{https://stacks.math.columbia.edu/tag/0539}{0539}}, of relative dimension $d$ by \cite{SP}*{Lemma~\href{https://stacks.math.columbia.edu/tag/0D4J}{0D4J}}, and even of pure relative dimension $d$ by \cite{SP}*{Lemma~\href{https://stacks.math.columbia.edu/tag/02FZ}{02FZ}}. In particular, by \cite{SP}*{Lemma \href{https://stacks.math.columbia.edu/tag/0AFE}{0AFE}} the local rings of $\ov{X}$ are all of dimension $\le d + 1$, so for an $x \in X$ of height $h$, every proper closed subset of the closure $\ov{\{x\}} \subset \ov{X}$ is of dimension $\le d - h$ (and is even of dimension $\le d - h - 1$ if $\mathcal O$ is $0$-dimensional). Since the closure $\ov{Z} \subset \ov{X}$ of $Z$ is the union of the closures of the generic points of $Z$, all of which are of height $\ge 2$ in $X$, this means that $\ov{Z} \setminus Z$ is $\cO$-fiberwise of codimension $\ge 2$  in $\ov{X}$ (and is even of codimension $\ge 3$ if $\mathcal O$ is $0$-dimensional). Similarly, letting $\sY$ be the closure of $Y_\fm$ in $\ov{X}_\fm$ where $\fm \subset \Spec \cO$ is the union of the closed points, we find that $\sY \setminus Y_\fm$ is $\mathcal O$-fiberwise of codimension $\ge 2$ in $X$. We replace the very ample line bundle $\cO_{\ov{X}}(1)$ by its large power and apply \cite{EGAIII1}*{corollaire~2.2.4} to force each global section of $\cO_{\ov{X}_\fm}(n)$ to lift to a global section of $\cO_{\ov{X}}(n)$ for $n > 0$. By applying \cite{split-unramified}*{Proposition 3.6} (especially, its last aspect to handle disconnected $\fm$; the $W$ there is our $X$ and the $Y$ there is our $\overline{Z}_\fm \cup \sY$) to the closed $\cO$-fibers of $\ov{X}$ and lifting the sections obtained to $\ov{X}$, we may even choose this large power so that there exist nonzero 
\[
h_0 \in \Gamma(\ov{X}, \sO_{\ov{X}}(1)), \ \ h_1 \in \Gamma(\ov{X}, \sO_{\ov{X}}(w_1)), \ \  \dotsc, \ \  h_{d - 1} \in \Gamma(\ov{X}, \sO_{\ov{X}}(w_{d - 1})) \qxq{with} w_1, \dotsc, w_{d - 1} > 0
\]
such that the hypersurfaces $H_i \ce V(h_i) \subset \ov{X}$ satisfy the following properties. 
\benumr
\m \label{m:PL-i}
$H_0$ does not contain $x_1, \dotsc x_n$.

\m \label{m:PL-ii}
The map $\pi\colon \ov{X} \setminus H_0 \ra \bA^{d - 1}_\cO$ determined by the $h_1/h_0^{w_1}, \dotsc, h_{d - 1}/h_0^{w_{d-1}}$ is smooth at each $x_i$.

\m\label{m:PL-iii}
$(\ov{Z} \cup \sY) \cap H_0 \cap \dotsc \cap H_{d - 1} = \emptyset$, in other words, $\ov{Z} \cup \sY$ does not meet the exceptional locus of the weighted blowing up in the following diagram determined by the $h_0, \dotsc, h_{d - 1}$:
\[
\xymatrix{
\ov{X} \setminus H_0 \ar[d]_-{\pi} \ar@{^(->}[r] &\Bl_{\ov{X}}(h_0, \dotsc, h_{d - 1}) \ar[d]_-{\ov{\pi}} \\
\bA^{d - 1}_{\cO} \ar@{^(->}[r] & \bP_\cO(1, w_1, \dotsc, w_{d - 1})
}
\]
(see \cite{split-unramified}*{Section 3.5} for a review of the weighted blowup $\Bl_{\ov{X}}(h_0, \dotsc, h_{d - 1})$; its formation may not commute with base change to $\fm$, but the formation of $\pi$ does).

\m \label{m:PL-extra}
Each $(\ov{Z} \cup \sY) \cap \ov{\pi}\i(\pi(x_i))$ lies in $\ov{X} \setminus H_0$.

\m\label{m:PL-iv}
In fact, each $(\ov{Z} \cup \sY) \cap \ov{\pi}\i(\pi(x_i))$ also lies both in $X$ and in the smooth locus of $\pi$. 
\eenum
By \ref{m:PL-iii}, each $(\ov{Z} \cup \sY) \cap \ov{\pi}\i(\pi(x_i))$ is a projective subscheme of $\ov{X}$, in fact, by \ref{m:PL-extra}, it is even a finite collection of possibly nonreduced points: indeed, any component of dimension $> 0$ would still be projective, and so could not lie in $\ov{X} \setminus H_0$ because the latter is affine. Thus, since $\ov{\pi}$ is projective, by spreading out and the openness of the quasi-finite locus of a morphism \cite{SP}*{Lemma~\href{https://stacks.math.columbia.edu/tag/01TI}{01TI}} applied to the projective morphism $\ov{\pi}|_{\ov{Z} \cup \sY}$, there is an affine open $S \subset \bA^{d - 1}_\cO$ containing all the $\pi(x_i)$ such that $(\ov{Z} \cup \sY) \cap \ov{\pi}\i(S)$ is $S$-quasi-finite, and hence, being projective, is even $S$-finite. By \ref{m:PL-extra}, at the cost of shrinking $S$ around the $\pi(x_i)$, we may then also ensure that $(\ov{Z} \cup \sY) \cap \ov{\pi}\i(S) = (\ov{Z} \cup \sY) \cap \pi\i(S)$. At the cost of further shrinking $S$ around the $\pi(x_i)$, we may then choose an affine open $X' \subset X \cap \pi\i(S)$ in the smoothness locus of $\pi$ containing all the $x_i$ and all the $(\ov{Z} \cup \sY) \cap \ov{\pi}\i(\pi(x_i))$ to make sure that even $(\ov{Z} \cup \sY) \cap X'$ is $S$-finite (it suffices to first choose any affine open $X'$ containing the indicated points and then base change to an affine open of $S$ containing all the $\pi(x_i)$ and not meeting the image of $((\ov{Z} \cup \sY) \cap \pi\i(S)) \setminus X'$, noting that this image is automatically closed due to finiteness).

Since $(\ov{Z} \cup \sY) \cap X' = (Z \cup Y_\fm) \cap X'$, we get that $Z \cap X'$ is also  $S$-finite. Thanks to \cite{SP}*{Lemma~\href{https://stacks.math.columbia.edu/tag/01TI}{01TI}} again, we may then shrink $S$ around the $\pi(x_i)$ and replace $X'$ by a suitable affine open containing all the $x_i$ and all the $(Z \cup Y_\fm) \cap \ov{\pi}\i(\pi(x_i))$ to also make $Y \cap X'$ be $S$-quasi-finite (in addition to $Z \cap X'$ being $S$-finite, as ensured by repeating the parenthetical argument at the end of previous paragraph). It remains to note that our smooth map $X' \ra S$ is of relative dimension $1$ by a dimension count. 
\epf

The following reembedding lemmas will help us to pass from the relative curve $X' \ra S$ of \Cref{lem:pres-lem} to a relative affine line. They are more subtle than the versions given in \cite{split-unramified}*{Lemma~6.3} or in prior references that developed the geometric approach to the Grothendieck--Serre conjecture because now $Y$ is merely quasi-finite. Relatedly, we do not know how to arrange that  $V = \bA^1_A$.

\blem \label{lem:embed-Y}
Let $Y$ be a quasi-finite, separated scheme over  a semilocal ring $A$ and, for each maximal ideal $\fm \subset A$, let $\iota_\fm \colon Y_{k_\fm} \hra \bA^1_{k_\fm}$ be a closed $k_\fm$-immersion. There are  principal affine opens $Y' \subset Y$ and $V \subset \bA^1_A$, both containing all the $Y_{k_\fm}$, and a closed immersion $\iota \colon Y' \hra V$ extending all the $\iota_\fm$. 
\elem

\bpf
Zariski Main Theorem \cite{EGAIV4}*{corollaire 18.12.13} gives an open immersion $Y \hra \wt{Y}$ into an $A$-finite scheme $\wt{Y} = \Spec(\wt{A})$. The union of the $Y_{k_\fm}$ is a closed subscheme of $\wt{Y}$ disjoint from $\wt{Y} \setminus Y$. Thus, some $a_\infty \in \wt{A}$ vanishes on $\wt{Y} \setminus Y$ and is a unit on every $Y_{k_\fm}$, and some $a \in \wt{A}$ is a unit on $\wt{Y} \setminus Y$ and is such that $a/a_\infty$ on each $Y_{k_\fm}$ is the $\iota_\fm$-pullback of the standard coordinate of $\bA^1_{k_\fm}$. Jointly, $a, a_\infty$ do not vanish at any point of $\wt{Y}$, so they determine a map $\wt{\iota} \colon \wt{Y} \ra \bP^1_A$ such that $\{a_\infty = 0\}$ set-theoretically is the $\wt{\iota}$-pullback of infinity. By construction, $\wt{\iota}$ extends the $\iota_\fm$ and $\wt{\iota}^{\,-1}(\bA^1_A) \subset Y$.   

The schematic image of $\wt{\iota}$ is an $A$-finite closed subscheme $\ov{Y} \subset \bP^1_A$: this is simpler when $A$ is Noetherian, but in general and more concretely, $\wt{\iota}$ factors through the affine complement $\Spec (B)$ of a hypersurface in $\bP^1_A$ disjoint from $\wt{\iota}(\wt{Y})$ (such a hypersurface exists by the avoidance lemma \cite{GLL15}*{Theorem~5.1}), $\overline Y$ is cut out by $\Ker(B \ra \wt A)$ because $\wt \iota$ is quasi-compact (see \cite{SP}*{Lemma~\href{https://stacks.math.columbia.edu/tag/01R8}{01R8}}), and the coordinate ring $\ov{A}$ of $\ov{Y}$ is the image of $B \rightarrow \wt{A}$, which is automatically $A$-finite because $B$ is of finite type over $A$ and $\wt{A}$ is $A$-finite. Thanks to this description, the image of the map $\wt{Y} \ra \ov{Y}$, which is finite by \cite{SP}*{Lemma \href{https://stacks.math.columbia.edu/tag/035D}{035D}}, contains every minimal prime of $\ov{A}$, so this map is surjective. In particular, for every maximal ideal $\fm \subset A$, the intersection $\ov{Y} \cap \bA^1_{k_\fm}$ set-theoretically is $Y_{k_\fm}$ (viewed inside $\bA^1_{k_\fm}$ via the monomorphism $\iota_\fm$), to the effect that, by the construction of $\wt \iota$, the finite map
\be \label{eqn:iota-tilde}
\wt{\iota}^{\,-1}(\bA^1_A) \ra \ov{Y} \cap \bA^1_A
\ee
is a closed immersion on $k_\fm$-fibers. By the Nakayama lemma \cite{SP}*{Lemma \href{https://stacks.math.columbia.edu/tag/00DV}{00DV} (6)}, this finite surjection that is injective on coordinate rings becomes also surjective on coordinate rings after semilocalizing $\ov{Y} \cap \bA^1_A$ along the union of its $k_\fm$-fibers. Thus, \eqref{eqn:iota-tilde} becomes an isomorphism after this semilocalization, so, by a limit argument, there is a principal affine open of $\ov{Y} \cap \bA^1_A$ containing its $k_\fm$-fibers over which the map \eqref{eqn:iota-tilde} is an isomorphism. This means that, as claimed, there are a principal affine open $Y' \subset \wt{\iota}^{\,-1}(\bA^1_A) \subset Y$ containing all the $Y_{k_\fm}$, a principal affine open $V \subset \bA^1_A$, and a closed immersion $\iota \ce \wt{\iota}|_{Y'} \colon Y' \hra V$ extending the $\iota_\fm$.
\epf

To use \Cref{lem:embed-Y} in practice, we need a criterion for the existence of the closed immersions $\iota_\fm$. \Cref{lem:embed-field} below gives such a criterion in terms of the following set-theoretic obstruction.

\bd
Let $A$ be a ring, let $Y$ be a quasi-finite $A$-scheme, and let $X$ be an $A$-scheme. There is no \emph{finite field obstruction} to embedding $Y$ into $X$ if for every maximal ideal $\fm \subset A$ with $k_\fm$ finite and every finite field extension $k'/k_\fm$, the number of $k'$-points of $Y_{k_\fm}$ does not exceed that of $X_{k_\fm}$.
\ed

The condition is fibral, but it is convenient to allow an arbitrary $A$ to simply be able to say that there is no finite field obstruction to embedding $Y$ into $X$ over $A$.

\blem \label{lem:embed-field}
For a finite scheme $Y$ over a field $k$ and a nonempty open $V \subset \bA^1_k$, there is a closed $k$-immersion $\iota\colon Y \hra V$ iff there is no finite field obstruction to it and $Y$ is a closed subscheme of \emph{some} smooth $k$-curve $C$, in which case we may choose $\iota$ to extend any $\iota_0 \colon Y_0 \hra V$ for a closed~$Y_0 \subset Y$. 
\elem

\bpf
The `only if' is clear, so we fix closed immersions $Y \subset C$ and $\iota_{0}$ as in the statement and assume that there is no finite field obstruction. We may build $\iota$ one connected component of $Y$ at a time and shrink $V$ at each step, so we may assume that $Y$ is connected with unique closed point $y$. If $k$ is finite, then the absence of the finite field obstruction allows us to choose a closed immersion $\iota_y\colon y \hra V$. If $k$ is infinite, then, since every closed point of a smooth curve over $k$ is also a closed point of $\bA^1_k$ (see \cite{split-unramified}*{Lemma 6.2}), we have a closed immersion $\iota_y\colon y \hra \bA^1_k$, and the possibility to change coordinates via $t \mapsto t + \gA$ for $\gA \in k$ allows us to assume that $\iota_y$ factors through $V$. In other words, for all $k$ we have reduced to the case when $Y_0 \neq \emptyset$.

In the case when the extension $k_y/k$ is separable, \cite{EGAIV4}*{proposition~17.5.3} 
ensures that the $n$-th infinitesimal neighborhood of $y$ in $C$ is $k$-isomorphic to $Y_n \ce \Spec(k_y[x]/(x^{n + 1}))$ over $k$ (the separability ensures that $k_y \tensor_k k_y$ has $k_y$ as a direct factor, so, by the invariance of the \'etale site under nilpotents, it suffices to identify the $n$-th infinitesimal neighborhood after base changing $C$ along $k \ra k_y$, that is, after reducing to the case $k \cong k_y$, in which \emph{loc.~cit.}~applies). This does not depend on $C$, and $Y \simeq Y_n$ for some $n \ge 0$. Thus, to extend our fixed $\iota_{0}$ to a desired $\iota$, by induction on $n \ge 0$, we only need to argue that every $k$-automorphism of $Y_m$ lifts to a $k$-automorphism of $Y_{m + 1}$. For this, by base change along the inverse of the induced $k$-automorphism of $k_y$, we may reduce to the case when this induced automorphism is the identity of $k_y$. This makes the automorphism $k_y$-linear, so we may replace $k$ by $k_y$ and further reduce to the case when $k_y= k$. In this case, however,  $k$-automorphisms of $Y_m$ correspond to elements $a_1x + \dotsc + a_mx^m \in k[x]/(x^{m + 1})$ with $a_i \in k$ and $a_1 \neq 0$, and such elements~lift.



In the remaining case when $k_y$ (equivalently, $k$) is infinite and $Y_0 \neq \emptyset$, it suffices to show that a given closed immersion $\iota_0 \colon Y_0 \hra V$ extends to a closed immersion of the square-zero infinitesimal neighborhood $\eps_{Y_0}$ of $Y_0$ in $C$: by iterating this with $Y_0$ replaced by $\eps_{Y_0}$ and eventually restricting to $Y$, we will obtain the desired $\iota$. By deformation theory, more precisely, by \cite{Ill05}*{Theorem 8.5.9~(a)}, the $k$-morphisms $\eps_{Y_0} \ra V$ that restrict to $\iota_{0}$ are parametrized by some affine space $\bA^N_k$. Since $\eps_{Y_0}$ is $k$-finite, the Nakayama lemma \cite{SP}*{Lemma \href{https://stacks.math.columbia.edu/tag/00DV}{00DV}} ensures that the locus parametrizing those $\eps_{Y_0} \ra V$ that are closed immersions is an open $\sV \subset \bA^N_k$. Moreover, $\sV \neq \emptyset$: indeed, we may check this after base change to any field extension of $k$, and a suitable such base change reduces us to the already settled case when $k_y/k$ is separable. Since $k$ is infinite and $\sV \subset \bA^N_k$ is nonempty, $\sV(k) \neq \emptyset$. Any $k$-point of $\sV$ corresponds to a sought closed immersion $\eps_{Y_0} \hra V$ that restricts to $\iota_0$. 
\epf

The embedding lemmas above help us build the following excision squares that allow us to pass~to~$\bA^1_A$.

\blem \label{lem:reembed}
Let $C$ be a smooth, affine scheme of pure relative dimension $1$ over a semilocal ring $A$, let $Y \subset C$ be an $A$-quasi-finite closed subscheme, and let $\iota_\fm \colon Y_{k_\fm} \hra \bA^1_{k_\fm}$ for maximal ideals $\fm \subset A$ be closed immersions. There are an affine open $C' \subset C$ containing the $Y_{k_\fm}$, an affine open $V \subset \bA^1_A$, and an \'{e}tale $A$-morphism $f\colon C' \ra V$ that embeds $Y \cap C'$ as a closed $Y' \subset V$ in such a way that
\[
\xymatrix{
Y \cap C' \ar@{^(->}[r] \ar[d]^-{\sim} & C' \ar[d]^-f \\
Y' \ar@{^(->}[r] & V
}
\]
is a Cartesian square in which the left vertical arrow is an isomorphism, as indicated. 
\elem
 
\bpf
By the final aspect of \Cref{lem:embed-field}, any fixed $\iota_\fm$ may be extended to any infinitesimal thickening of $Y_{k_\fm}$ in $C_{k_\fm}$. In particular, we lose no generality by replacing $Y$ by any of its infinitesimal neighborhoods in $C$, so we may and do assume that each clopen of every $Y_{k_\fm}$ is nonreduced. By \Cref{lem:embed-Y}, there are principal affine opens $Y' \subset Y$ and $V \subset \bA^1_A$, both containing all the $Y_{k_\fm}$, and a closed immersion $\iota \colon Y' \hra V$ extending the $\iota_\fm$. Since $Y' \subset Y$ is a principal affine open, we may replace $C$ by a principal affine open containing all the $Y_{k_\fm}$ to reduce to $Y' = Y$. By lifting the $\iota$-pullback of the standard coordinate of $\bA^1_A$, we then extend $\iota\colon Y \hra V$ to an $A$-morphism~$f \colon C \ra \bA^1_A$. 

By \cite{SP}*{Lemma \href{https://stacks.math.columbia.edu/tag/01TI}{01TI}}, the quasi-finite locus of $f$ is open, and the $A$-smoothness of $C$ together with the nonreducedness of each clopen of every $Y_{k_\fm}$ force this locus to contain all the $Y_{k_\fm}$: indeed, if $C'$ is an irreducible component of $C_{k_\fm}$ containing  a point of $Y_{k_\fm}$, then $f|_{C'}$ is quasi-finite because $f$ cannot collapse $C'$ to a point of $\bA^1_{k_\fm}$ since $f|_{Y_{k_\fm}}$  is a closed immersion and the components of $Y_{k_\fm}$ are nonreduced. Moreover, since $C$ and $\bA^1_A$ are $A$-smooth, we may $A$-fiberwise apply the flatness criterion \cite{EGAIV2}*{proposition 6.1.5} to see that $f$ is $A$-fiberwise flat at the points of each $Y_{k_\fm}$. Thus, \cite{EGAIV3}*{corollaire~11.3.11} implies that $f$ itself is flat at the points of each $Y_{k_\fm}$. Since \'{e}taleness of a flat morphism may be checked fiberwise and all the components of all the $Y_{k_\fm}$ are nonreduced, $f|_{Y_{k_\fm}}$ being closed immersions then implies that all the $Y_{k_\fm}$ even lie in the \'{e}tale locus of $f$. Consequently, we may replace $Y$ and $V$, and then also $C$, by principal affine opens containing all the $Y_{k_\fm}$ to reduce to the case when $f$ is \'{e}tale. 

Finally, we replace $C$ by the $f$-preimage of $V$ to make $f$ factor through $V$. Since a section of separated \'{e}tale morphism is a clopen immersion, $f\i(f(Y)) = Y \sqcup \wt{Y}$ for some closed $\wt{Y} \subset C$. By inverting a function on $C$ that vanishes on $\wt{Y}$ but is a unit on $Y$, we get a desired affine open $C'$. 
\epf

We are ready to build the promised closed $Z \subset \Spec R$ of codimension $\ge 2$ away from which our $G$-torsor is simpler: any $V$ as in the following proposition contains $\bP^1_{\Spec(R) \setminus Z}$ for some such $Z$. 

\bprop \label{prop:build-Z}
Let $G$ be a reductive group over a Noetherian semilocal ring  $R$ that is flat and geometrically regular over some Dedekind ring. For a generically trivial $G$-torsor $E$ over $R$, 
there are
\benumr

\m \label{m:BZ-ii}
an open $V \subset \bP^1_R$ containing all the height $\le 2$ points and the sections $\{t = 0\}$ and $\{t = \infty\}$\uscolon

\m \label{m:BZ-iii}
a $G$-torsor $\sE$ over $V$ that trivializes away from some $R$-quasi-finite closed of $V$ and is such~that 
\[
\qq \sE|_{\{t = 0\}} \simeq E,  \qxq{$\sE|_{\{t = \infty\}}$ is trivial, and} \x{$\sE|_{\bP^1_{\Frac(R)}}$ is trivial}.
\]
\eenum
\eprop

\bpf
We first dispose of the condition that $V$ cover the height $\le 2$ points, so we suppose that $V \subset \bP^1_R$ is an open satisfying the other conditions, in particular, such that $\sE$ trivializes away from an $R$-quasi-finite closed $\sY \subset V$ and also on $V_{\Frac(R)}$. By spreading out, $\sE$ is trivial over $V_S$ for some dense open $S \subset \Spec(R)$. By patching with the trivial torsor over $\bP^1_S$, we may assume that $V$ contains $\bP^1_S$, so also contains $\bP^1_{\Frac(R)}$ and is trivial thereon. The closure $\ov{\sY} \subset \bP^1_R$  of $\sY$ is $R$-finite because $V \subset \bP^1_R$ is $R$-fiberwise dense. By patching, $\sE$ over $V$ extends to a $G$-torsor over $V \cup (\bP^1_R \setminus \ov{\sY})$ that trivializes away from the $R$-quasi-finite closed $\sY$. Thus, we replace $V$ by $V \cup (\bP^1_R \setminus \ov{\sY}) = \bP^1_R \setminus (\ov\sY \setminus \sY)$ to force $V$ to cover the height $\le 1$ points of $\bP^1_R$.  At this point, by \cite{CTS79}*{th\'{e}or\`{e}me~6.13} (see also \cite{torsors-regular}*{Section~1.3.9} for a recap) applied to the local rings of the generic points of $\bP^1_R \setminus V$ and then spreading out and patching, $\sE$ over $V$ extends to a $G$-torsor over some open of $\bP^1_R$ covering the height $\le 2$ points. Consequently, we may enlarge $V$ again to cover the height $\le 2$ points.

Having disposed of the codimension requirement, we let $\cO$ be a Dedekind ring over which $R$ is flat and geometrically regular and decompose $\cO$ and $R$ into factors to force both of them to be domains. We then combine Popescu's \cite{SP}*{Theorem~\href{https://stacks.math.columbia.edu/tag/07GC}{07GC}} with a limit argument to reduce to the case when $R$ is the semilocal ring of a smooth, affine, integral $\cO$-scheme $X$. We spread out to assume that $G$ and $E$ begin life over $X$. We may assume that $X$ is of relative dimension $d > 0$ over $\cO$ because else $E$ is trivial by the settled Dedekind case of the Grothendieck--Serre conjecture, see \cite{Guo22a}*{Theorem~1}, a case in which we may choose $V = \bP^1_R$ with $\sE$ trivial. 

More generally, we apply \cite{Guo22a}*{Theorem 1} to the semilocalization of $X$ at the union of the generic points of the closed $\cO$-fibers of $X$ and use a limit argument to find a closed $Y \subset X$ that contains no irreducible component of any $\cO$-fiber of $X$ and is such that $E$ trivializes over $X \setminus Y$. By \Cref{lem:pres-lem} (applied with $Z = \emptyset$) and the fact that every open immersion $S \subset \mathbb A^{d - 1}_\cO$ is quasi-finite, at the cost of shrinking $X$ around $\Spec R$ we may find a smooth morphism $X \ra \bA^{d - 1}_\cO$ of relative dimension $1$ with respect to which $Y$ is quasi-finite. Base change along $\Spec R \ra \bA^{d - 1}_\cO$~then~gives
\bitem
\m
a smooth, affine $R$-scheme $C$ of pure relative dimension $1$ equipped with an $s \in C(R)$\uscolon 

\m
a reductive $C$-group scheme $\sG$ with $s^*(\sG) \cong G$ and a $\sG$-torsor $\sE$ over $C$ with $s^*(\sE) \cong E$\uscolon

\m
an $R$-quasi-finite closed subscheme $\sY \subset C$ containing $s$ such that $\sE|_{C \setminus \sY}$ is trivial.
\eitem
We will gradually simplify the data of these $C$, $s$, $\sG$, $\sE$, and $\sY$ over $R$ to arrive at our $V \subset \bP^1_R$. 

By \cite{Li24}*{Proposition 7.4} and spreading out, there are a finite \'{e}tale cover $\wt{C} \surjects C'$ of some affine open neighborhood $C' \subset C$ of $s$, a lift $\wt{s} \in \wt{C}(R)$ of $s$, and a reductive group isomorphism $\sG_{\wt{C}} \cong G_{\wt{C}}$ whose $\wt{s}$-pullback agrees with the identification $s^*(\sG) \cong G$. By replacing $(C, s)$ by $(\wt{C}, \wt{s})$ and $\sG$, $\sE$, $\sY$ by their pullbacks to $\wt{C}$, we therefore reduce to the case when $\sG \cong G_{\wt{C}}$. 

Since $\sY$ is merely required to be $R$-quasi-finite (and not $R$-finite), we may replace $C$ by some affine open containing $s$ to arrange that set-theoretically $\sY_{k_\fm} = s_{k_\fm}$ for every maximal ideal $\fm \subset R$. Since $\#\bA^1_R(k_\fm) \ge 2$, this ensures that there is no finite field obstruction to embedding $\sY \sqcup \Spec R$ into $\bA^1_R$. Therefore, \Cref{lem:embed-field,lem:reembed} give us an affine open $C' \subset C \sqcup \bA^1_R$ containing $s \sqcup \{t = 0\}$, an affine open $V \subset \bA^1_R$, and an \'{e}tale $R$-morphism $f \colon C' \ra V$ that fits into a Cartesian square
\[
\xymatrix{
(\sY \cap C') \sqcup \{ t = 0 \} \ar@{^(->}[r] \ar[d]^-{\sim} & C' \ar[d]^-{f} \\
\sY' \ar@{^(->}[r] & V
}
\]
for some closed subscheme $\sY' \subset V$. By patching the disjoint union of $\sE$ over $C' \cap C$ and the trivial $G$-torsor over $C' \cap \bA^1_R$ with the trivial $G$-torsor over $V \setminus \sY'$ (see, for instance, \cite{torsors-regular}*{Proposition~4.2.1}), we therefore obtain a $G$-torsor $\sE'$ over $V$ such that $\sE'|_{V \setminus \sY'}$ is trivial and disjoint $s, s_0 \in V(R)$ such that $s^*(\sE') \cong E$ and $s_0^*(\sE')$ is trivial. By \cite{Gil02}*{corollaire 3.10 (a)}, the triviality away from $\sY'$ implies that $\sE'$ is also trivial over $V_{\Frac(R)}$. 

At this point, in the view of the initial reduction described in the first paragraph of the proof, we have basically already constructed all the required data. To finish, we note that since $R$ is semilocal, the automorphism group of $\bP^1_R$ acts transitively on $\bP^1_R(R)$. Thus, we may assume that $s_0$ is the $R$-point $\{t = \infty\}$. Since $s_0$ is disjoint from $s$, we then shift the coordinate of $\bA^1_R$ to arrange that, in addition, $s$ is the $R$-point $\{t = 0\}$.
\epf


\section{Torsors over $\bP^1_A$ via the geometry of $\Bun_G$} \label{sec:BunG}

To proceed further, we need to analyze the $G$-torsor $\sE$ over $V \subset \bP^1_R$ obtained in \Cref{prop:build-Z}. An initial step to this and a general bedrock of the geometric approach to the Grothendieck--Serre conjecture is the fact that a $G$-torsor on $\bP^1_A$ over a semilocal ring is $A$-sectionwise constant. This constancy was recently established by Panin--Stavrova in \cite{PS25}, and we reprove and mildly generalize their result in \Cref{thm:section-const} below. The constancy comes from the following geometric property of the algebraic stack $\Bun_G$ parametrizing $G$-bundles on $\bP^1_A$, in addition, \Cref{prop:const} simultaneously reproves, strengthens, and explains its numerous special cases in \cite{PSV15}*{Proposition~9.6}, \cite{Tsy19}, \cite{Fed21a}*{Proposition~2.2}, \cite{split-unramified}*{Lemma~8.3}, and elsewhere. For a basic review of some properties of algebraic stacks that are useful for studying torsors, see \cite{topology-torsors}*{Appendix A}.
 
\bprop \label{prop:const}
Let $\pi \colon C \ra S$ be a proper, flat, finitely presented scheme morphism and let $G$ be a flat, finitely presented, quasi-affine $S$-group. The restriction of scalars $\Bun_G \ce \pi_*((\bbB G)_C)$ is a locally finitely presented algebraic $S$-stack with quasi-affine diagonal. The adjunction morphism
	\[
\bbB G \ra \Bun_G
\]
\benum
\m \label{m:const-a}
is a monomorphism of algebraic $S$-stacks if $H^0(C_s, \sO_{C_s}) \cong k_s$ for $s \in S$\uscolon 

\m \label{m:const-b}
is an open immersion if $H^0(C_s, \sO_{C_s}) \cong k_s$ and $H^1(C_s, \sO_{C_s}) \cong H^2(C_s, \sO_{C_s}) \cong 0$ for $s \in S$.
\eenum
When \ref{m:const-b} holds with $S$ quasi-compact, a $G$-torsor over $C$ descends to $S$ iff it does so on the closed $S$-fibers~of~$C$. 

\eprop

The main case of interest for us is $C = \bP^1_S$ but the proof is no more difficult in general. 

\bpf
By \cite{SGA3Inew}*{expos\'e VI$_B$, proposition 11.11 (i)$\Leftrightarrow$(ii)}, the quasi-affine $S$-group $G$ has affine fibers, so the algebraic stack $\bbB G$ has affine stabilizers. Thus, since $\bbB G$ is finitely presented and has a quasi-affine diagonal (see \cite{topology-torsors}*{Lemma A.2 (b)}), the representability of $\Bun_G$ by an algebraic stack, as well as its geometric properties, follow from \cite{HR19}*{Theorem 1.3}. Moreover, the final claim about the closed $S$-fibers follows from \ref{m:const-b} because any open containing all the closed points of a quasi-compact scheme is the entire scheme (equivalently, a quasi-compact scheme has a closed~point).

 In \ref{m:const-a}, by base change and \cite{SP}*{Lemma~\href{https://stacks.math.columbia.edu/tag/04ZZ}{04ZZ}}, it suffices to check the full faithfulness of $\bbB G \ra \Bun_G$ on $S$-points. For this, for any  $G$-torsors $E$ and $E'$ over $S$, we need to check that 
\be \label{eqn:isom-equal}
\underline{\Isom}_G(E, E')(S) \isomto \underline{\Isom}_G(E, E')(C).
\ee
By working fpqc locally on $S$ to trivialize $E$ and $E'$, it is enough to argue that $G(S) \isomto G(C)$ and, by also using \cite{EGAIV4}*{corollaire 17.16.2}, we may assume that $C(S) \neq \emptyset$, so that $G(S) \hra G(C)$. For the surjectivity, we may again work locally and now combine Noetherian approximation (with \cite{Ill05}*{Corollary 8.3.11 (a)} to keep the assumption on $H^0$) with the rigidity lemma \cite{MFK94}*{Proposition~6.1} 
to reduce to the case when $S$ is the spectrum of a field $k$. In the field case, however, since morphisms to an affine scheme correspond to ring homomorphisms induced on global sections,  the assumption $H^0(C, \sO_C) \cong k$ and the affineness of $G$ imply that every $C$-point of $G$ descends to a $k$-point. 

In \ref{m:const-b}, we already know from \ref{m:const-a} that the map is a monomorphism, and hence is representable by algebraic spaces by \cite{SP}*{Lemmas \href{https://stacks.math.columbia.edu/tag/04Y5}{04Y5} and \href{https://stacks.math.columbia.edu/tag/04ZZ}{04ZZ}}. Thus, it suffices to check that it is formally smooth: indeed, it will then be smooth by \cite{SP}*{Lemmas \href{https://stacks.math.columbia.edu/tag/06Q6}{06Q6} and \href{https://stacks.math.columbia.edu/tag/0DP0}{0DP0}}, hence representable by schemes by Rydh's \cite{SP}*{Lemmas \href{https://stacks.math.columbia.edu/tag/0B8A}{0B8A}}, and so an open immersion by \cite{SP}*{Theorem~\href{https://stacks.math.columbia.edu/tag/025G}{025G}}. Concretely, for the formal smoothness, given a square-zero thickening $T \hra T'$ of affine $S$-schemes, we need to argue that a $G$-torsor $\sE$ over $C_{T'}$ descends to $T'$ granted that its restriction to $C_T$ descends to a $G$-torsor $E$ over $T$. Let $J \subset \cO_{T'}$ be the ideal sheaf of $T$, so that $J^2 = 0$ and we may view $J$ as a quasi-coherent $\cO_T$-module. By \ref{m:const-a}, we already know that, if a sought descent exists, it is unique up to a unique isomorphism, so we may work fpqc locally on $T'$ to assume that 
\[
H^1(C_T, \sO_{C_T}) \cong H^2(C_T, \sO_{C_T})  \cong 0
\]
(see \cite{Ill05}*{Corollary 8.3.11}), that the co-Lie complex $\ell_{E/T}$, controlling the deformations of $E$, consists of free vector bundles placed in degrees $-1$ and $0$ (see \cite{Ill72}*{\'equation (2.4.2.9), page~208}), and, as in \ref{m:const-a}, that $C(T') \neq \emptyset$. By \cite{Ill05}*{equation (8.3.2.2) and Corollary 8.3.6.5 (a)} (we apply the corollary to $X \ce T$ and $E \ce R\Gamma(C_T, \sO_{C_T})$, with $M \ce J$), the displayed vanishing ensures that
\[
H^1(C_T, J|_{C_T}) \cong H^1(C_T, \sO_{C_T}) \tensor_{\sO_T} J \cong 0 \qxq{and}  H^2(C_T, J|_{C_T}) \cong H^2(C_T, \sO_{C_T}) \tensor_{\sO_T} J \cong 0.
\]
Consequently, the structure of $\ell_{E/T}$ forces the vanishing 
\[
\Ext^1_{\sO_{C_T}}(\ell_{E/T}|_{C_T}, J|_{C_T}) \cong 0.
\]
Thus, \cite{Ill72}*{th\'{e}or\`{e}me 2.4.4, page 209} implies that $\sE$ is the unique deformation of $E|_{C_T}$ to a $G$-torsor over $C_{T'}$. Since the pullback of $\sE$ along any $T'$-point of $C$ is another such deformation, $\sE$ must agree with this base change, so $\sE$ is constant. \epf


Even when $C = \bP^1_S$, the open immersion of \Cref{prop:const}~\ref{m:const-b} is typically not closed, for instance, this would contradict \cite{Fed16}*{Theorems 3~(ii) and 5}.  Nevertheless, it is closed when $G$ is of multiplicative type, as follows from the following broadly useful and widely known lemma that generalizes \cite{GR18}*{Proposition 11.4.2}, \cite{Fed22a}*{Lemma 2.14}, and other results in the literature.

\blem \label{lem:BunM}
For a finite type group $M$ of multiplicative type over a scheme $S$, its cocharacter $S$-scheme $X_*(M) \ce \underline{\Hom}_{\gp}(\bG_m, M)$, and the $S$-stack $\Bun_M$ parametrizing $M$-torsors over $\bP^d_S$ with $d > 0$, we have
\[
\Bun_M \cong \bbB M \times_S X_*(M), \qxq{in particular,} H^1(\bP^d_S, M) \cong H^1(S, M) \oplus H^0(S, X_*(M));
\]
if $M$ is, in addition, finite, then $\Bun_M \cong \bbB M$ and, in particular, $\bbB M(S) \isomto (\bbB M)(\bP^1_S)$.
\elem

 \bpf
 For finite $M$, we have $X_*(M) \cong 0$, so the claims about finite $M$ follow from the rest. 
 
The map $\bbB M \times_S X_*(M) \ra \Bun_M$ is given on $S$-points as follows: a pair of an $M$-torsor $E$ over $S$ and an $S$-morphism $\gA \colon \bG_{m,\, S} \ra M_{S}$ is sent to the contracted product\footnote{Since $M$ is commutative, the contracted product of two $M$-torsors $E_1$ and $E_2$ may be defined simply as the inflation of the $(M \times M)$-torsor $E_1 \times E_2$ to an $M$-torsor along the multiplication map $M \times M \ra M$.} of $E|_{\bP^d_{S}}$ and the extension along $\gA|_{\bP^d_{S}}$ of the $\bG_m$-torsor corresponding to $\sO(1)$, and similarly for points valued in a variable $S$-scheme $S'$. By the flexibility of base change to $S'$, it suffices to show that every $M$-torsor $\sE$ over $\bP^d_S$ arises from $E$ and $\gA$ as above that are uniquely determined up to a unique isomorphism.


Certainly, $E$ is uniquely determined by $E \simeq p^*(\sE)$ for a fixed $p \in \bP^d_S(S)$, so, by twisting and using the bijection $M(S) \isomto M(\bP^d_{S})$ that results as in \eqref{eqn:isom-equal}, all we need to show is that $\sE$ comes from a unique $\gA$ when $p^*(\sE)$ is trivialized. Due to this rigidification along $p$ and the fact that, by $M(S) \isomto M(\bP^d_{S})$, isomorphisms of rigidified $M$-torsors over $\bP^d_S$ are unique if they exist, the claim is fpqc local over $S$. Thus, we assume that $S = \Spec A$ is affine, then, by a limit argument, that $A$ is local, and, by decomposing $M$, that $M$ is either $\bG_{m,\, S}$ or $\mu_{n,\, S}$. For $\bG_m$, the desired $H^1(\bP^d_A, \bG_m) \cong \bZ$ holds when $A$ is a field, so, by \Cref{prop:const}, also when $A$ is local. The $\mu_n$ case follows from this by the sequence $0 \ra \mu_n \ra \bG_m \xra{n} \bG_m \ra 0$ and the isomorphism $\bG_m(A) \isomto \bG_m(\bP^d_A)$.
\epf

For finite groups $M$ of multiplicative type, we may slightly extend \Cref{lem:BunM} to gerbes as follows. We recall that an $M$-gerbe is a stack that fppf locally on the base is isomorphic to the stack $\bbB M$ of $M$-torsors and that up to equivalence $M$-gerbes are classified by $H^2_\fppf$ with coefficients in $M$, see \cite{Gir71}*{chapitre III, d\'efinition 2.1.1, section 2.1.1.2, corollaire 2.2.6; chapitre IV, th\'eor\`eme 3.4.2 (i)}.

\blem \label{lem:gerbeM}
Let $M$ be a finite group of multiplicative type over a scheme $S$ and fix a $d > 0$. 
\benum
\m \label{m:gM-a}
For an $M$-gerbe $\sM$ over $\bP^d_S$, the $s \in S$ such that $\sM$ trivializes over $\bP^d_{\ov{k}_s}$ form a clopen~$S_\sM \subset S$.

\m \label{m:gM-b}
Base change is an equivalence between the $(2, 1)$-category of $M$-gerbes over $S$ and that of those $M$-gerbes $\sM$ over $\bP^d_S$ with $S_\sM = S$\uscolon in particular, each $\sM$ trivializes fppf locally on~$S_\sM$. 
 
\eenum
\elem

\bpf
By descent, for both claims we may work fppf locally on $S$, so we may assume that $M$ is a product of various $\mu_{n,\, S}$, in particular, that there are split $S$-tori $T$ and $T'$ and an exact sequence
\[
0 \ra M \ra T \ra T' \ra 0.
\]
By \cite{Gab81}*{Chapter II, Part 2, Theorem 2 on page 193}, each element of $H^2(\bP^d_S, T)_\tors$ descends to $H^2(S, T)$. Thus, by \Cref{lem:BunM}, in \ref{m:gM-a} we may fppf localize $S$ further to reduce to the case when the class of $\sM$ in $H^2(\bP^d_S, M)$ comes from an $S$-point of the constant $S$-scheme $X_*(T')/X_*(T)$. By \Cref{lem:BunM} again, the locus of $S$ over which this $S$-point is the zero section is the sought $S_\sM$. Moreover, we have simultaneously shown the last aspect of \ref{m:gM-b}: $\sM$ trivializes fppf locally on $S_\sM$. 

For \ref{m:gM-b}, we first note that for any $S$-scheme $S'$, the $S'$-endomorphisms of the trivial $M$-gerbe $\bbB M$ are given by the contracted products with $M$-torsors over $S'$ (the relevant $M$-torsor over $S^\prime$ is simply the image of the trivial $M$-torsor under the endomorphism in question), to the effect that all such endomorphisms are automorphisms and their groupoid is identified with $(\bbB M)(S')$. Thus, the full faithfulness in \ref{m:gM-b} follows from fppf descent and the equivalence $(\bbB M)(S') \isomto (\bbB M)(\bP^d_{S'})$ supplied by \Cref{lem:BunM}. The essential surjectivity then follows from descent and the already established last aspect of \ref{m:gM-b}. 
\epf

The following lemma is useful for lifting the structure group of a torsor over $\bP^1_S$ along an isogeny $\wt{G} \ra G$. It is, of course, possible to analyze the geometry of the map $\Bun_{\wt{G}} \ra \Bun_G$ more thoroughly  but we do not pursue this here in order to keep our focus on what is needed for \Cref{thm:section-const}.

\blem  \label{lem:Cartesian}
For an isogeny $\wt{G} \ra G$ of reductive $S$-groups, the image of the map $\Bun_{\wt{G}} \ra \Bun_G$ between algebraic $S$-stacks parametrizing torsors over $\bP^d_S$ with $d > 0$ is clopen. For any $p \in \bP^d_S(S)$, the following square is Cartesian\ucolon
\[
\q \xymatrix{
\Bun_{\wt{G}} \ar[d]_-{\wt{\sE} \mapsto p^*(\wt{\sE})} \ar[r] & \im(\Bun_{\wt{G}} \ra \Bun_G) \ar[d]^-{\sE \mapsto p^*(\sE)} \\
\mathbf B\wt{G} \ar[r] & \mathbf B G,
}
\]
in particular, a $G$-torsor $\sE$ over $\bP^d_S$ lifts to a $\wt{G}$-torsor $\wt{\sE}$ iff it does so both on geometric $S$-fibers and after pullback by the $S$-point $p$, in which case giving $\wt{\sE}$ amounts to giving $p^*(\wt{\sE})$. 
\elem
 
\bpf
Set $M \ce \Ker(\wt{G} \ra G)$. For a $G$-torsor $\sE$ over an $S$-scheme $S'$, the category that parametrizes its liftings to a $\wt{G}$-torsor over variable $S'$-schemes is an $M$-gerbe over $S'$ (see \cite{topology-torsors}*{Proposition~A.4~(d) and its proof}), in particular, $\sE$ lifts to a $\wt{G}$-torsor iff this $M$-gerbe is trivial. Consequently, \Cref{lem:gerbeM}~\ref{m:gM-a} implies that image of the map $\Bun_{\wt{G}} \ra \Bun_G$ is clopen, whereas \Cref{lem:gerbeM}~\ref{m:gM-b} implies that the depicted square is indeed Cartesian.  
\epf

We turn to the promised $A$-sectionwise constancy of $G$-torsors over $\bP^1_A$ for semilocal $A$. Our argument for it is similar to that of the case treated by Panin--Stavrova in \cite{PS25}, even if perhaps slicker thanks to the geometric machinery above. In turn, their argument is slicker but somewhat similar to Fedorov's \cite{Fed22a}*{Theorem 6} that was mildly generalized in \cite{torsors-regular}*{Proposition 5.3.6}. The general idea goes back at least to \cite{PSV15}, \cite{FP15}, and \cite{Fed16}.

\bthm \label{thm:section-const}
 For a reductive group $G$ over a semilocal ring $A$, every $G$-torsor $\sE$ over $\bP^d_A$ is $A$-sectionwise constant\ucolon up to isomorphism, the $G$-torsor $s^*(\sE)$ over $A$ does not depend on $s \in \bP^d_A(A)$. 
 \ethm
 
 \bpf
We may assume that $d > 0$. The projective $d$-space over any field then has at least three rational points (even $\mathbb P^1_{\bF_2}$ has three distinct rational points!), so $\bP^d_A$ has an $A$-point that is disjoint from any two fixed $A$-points. Moreover, two disjoint $A$-points lie on a uniquely determined $\bP^1_A \subset \bP^d_A$. Thus, overall we may assume that $d = 1$. Moreover, since $A$ is semilocal, for any $s \in \bP^1_A(A)$, there is an $s' \in \bA^1_A(A)$ disjoint from $s$. Thus, we may change coordinates to first make $s^\prime$ be $\{t = 0\}$ and then make $s$ be $\{t = \infty\}$, and hence reduce to showing that $\sE|_{\{t = \infty\}} \simeq \sE|_{\{t = 0\}}$. By then replacing $G$ by an inner twist, it even suffices to show that $\sE|_{\{t = 0\}}$ is trivial granted that so is $\sE|_{\{t = \infty\}}$. 

Let $\sF$ be the $\Corad(G)$-torsor over $\bP^1_A$ obtained by inflating $\sE$. \Cref{lem:BunM} ensures that $\sF|_{\{t = 0\}}$ is trivial and that $\sF$ comes from an element of $X_*(\Corad(G))(A)$. Thus, since $\sO(1)$ pulls back to $\sO(d)$ under the map $\varphi_d\colon \bP^1_A \ra \bP^1_A$ that raises the homogeneous coordinates to their $d$-th powers, by choosing $d$ to be the degree of the isogeny $\Rad(G) \ra \Corad(G)$ and replacing $\sE$ by $\varphi_d^*(\sE)$ we reduce to the case when $\sF$ lifts to a $\Rad(G)$-torsor over $\bP^1_A$ that comes from an element of $X_*(\Rad(G))(A)$, in particular, that is $A$-sectionwise trivial. By twisting $\sE$ by this $\Rad(G)$-torsor, we therefore reduce to the case when $\sF$ is trivial. This means that $\sE$ lifts to a $G^\der$-torsor over $\bP^1_A$, to the effect that we have reduced to the case when $G$ is semisimple. This reduction might force us to revert to showing that $\sE|_{\{t = \infty\}} \simeq \sE|_{\{t = 0\}}$ without knowing that $\sE|_{\{t = \infty\}}$ is trivial, but we may afterwards again replace $G$ by an inner twist as in the first paragraph of the proof to still arrange that $\sE|_{\{t = \infty\}}$ be trivial. 

Once $G$ is semisimple, we pullback by $\varphi_d$ again, with $d$ now being the degree of the isogeny $G^{\mathrm{sc}} \ra G$: by \cite{Gil02}*{th\'{e}or\`{e}me~3.8}, this has the advantage of ensuring that each $\sE|_{\bP^1_{\ov{k}_s}}$ for $s \in S$ now lifts to a $G^{\mathrm{sc}}$-torsor over $\bP^1_{\ov{k}_s}$. By \Cref{lem:Cartesian}, then $\sE$ itself lifts to a $G^{\mathrm{sc}}$-torsor over $\bP^1_A$ whose restriction to infinity is trivial, to the effect that we have reduced to the case when $G$ is semisimple, simply connected. Due to \cite{SGA3IIInew}*{expos\'{e} XXIV, section 5.3, propositions 5.10 (i), 8.4} (that is, the analogue of \eqref{eqn:tot-isot}), we may then even assume that $G$ is simple. 

At this point, we begin the remaining argument by settling the isotropic case in the following claim. 

\bcl \label{cl:isotropic}
Let $A$ be a semilocal ring, let $G$ be a simple, simply connected $A$-group that is isotropic in the sense that it has an $A$-fiberwise proper parabolic $A$-subgroup, and let $\sE$ be a $G$-torsor over $\bP^1_A$. If $\sE|_{\{t = \infty\}}$ is trivial, then $\sE|_{\bA^1_A}$ is also trivial, so that $\sE|_{\{t = 0\}}$ is trivial, too. 
\ecl

\bpf
The assumptions on $G$ ensure that the following map is surjective:
\be \label{eqn:tot-iso-surj}
\tst G(A\llp t\i \rrp)/G(A\llb t\i \rrb) \surjects \prod_\fm G(k_\fm\llp t\i \rrp)/G(k_\fm\llb t\i \rrb),
\ee
where $\fm$ ranges over the maximal ideals of $A$, see \cite{split-unramified}*{(2) in the proof of Proposition 8.4} (the essential input here is the Borel--Tits theorem \cite{Gil09}*{fait 4.3, lemme 4.5}; the displayed surjectivity is also very close to \cite{Fed16}*{Proposition 7.1} and, implicitly, it is an important part of \cite{FP15}). Thanks to our assumption that $\sE|_{\{t = \infty\}}$ is trivial, Henselian invariance \cite{Hitchin-torsors}*{Theorem 2.1.6} ensures that $\sE$ is also trivial over $A\llp t\i \rrp$. Now by patching for $G$-torsors \cite{Hitchin-torsors}*{Lemma 2.2.11~(b)} or \cite{Fed16}*{Proposition 4.4}, the surjectivity \eqref{eqn:tot-iso-surj} means that every $G$-torsor over $\bigsqcup_\fm \bP^1_{k_\fm}$ that is obtained by patching $\sE|_{\bigsqcup_\fm \bA^1_{k_\fm}}$ with the trivial $G$-torsor at infinity lifts to a $G$-torsor over $\bP^1_A$ obtained by patching $\sE|_{\bA^1_A}$ with the trivial $G$-torsor at infinity. However, $\sE|_{\bigsqcup_\fm \bA^1_{k_\fm}}$ is trivial by \cite{Gil02}*{lemme~3.12}, so we get that $\sE|_{\bA^1_A}$ extends to a $G$-torsor $\sE'$ over $\bP^1_A$ such that $\sE'|_{\bigsqcup_\fm \bP^1_{k_\fm}}$ and $\sE'|_{\{t = \infty\}}$ are both trivial. By \Cref{prop:const}, then $\sE'$ itself is trivial, so that $\sE|_{\bA^1_A}$ is trivial, too. 
\epf

In the remaining case when our simple, simply connected $A$-group $G$ is not isotropic, let us consider any $A$-(finite \'{e}tale) subscheme $Y = \Spec A' \subset \bG_{m,\, A}$ such that $G_Y$ is isotropic and for each maximal ideal $\fm \subset A$ with $G_{k_\fm}$ isotropic, $Y_{k_\fm}$ has two disjoint nonempty clopens of coprime degrees over $k_\fm$ (we will later build such a $Y$). We may apply the settled isotropic case after base change along $Y \ra \Spec A$, so, since $Y \subset \bA^1_A$ gives rise to a $Y$-point of $\bA^1_Y$, we see that $\sE|_Y$ is trivial. On the other hand, \eqref{eqn:tot-iso-surj} applied after such a base change gives
\be \label{eqn:ani-surj}
\tst G(A'\llp y \rrp)/G(A'\llb y \rrb) \surjects \prod_\fm G((A'\tensor k_\fm)\llp y \rrp)/G((A'\tensor k_\fm)\llb y \rrb),
\ee
where $\fm$ still ranges over the maximal ideals of $A$. Since our choice of $Y$ and \cite{Gil02}*{th\'{e}or\`{e}me 3.8} still ensure that $\sE|_{\bigsqcup_\fm (\bP^1_{k_\fm} \setminus Y_{k_\fm})}$ is trivial, analogously to the previous paragraph, this surjectivity implies that $\sE|_{\bP^1_A \setminus Y}$ extends to a $G$-torsor $\sE'$ over $\bP^1_A$ such that $\sE'|_{\bigsqcup_\fm \bP^1_{k_\fm}}$ is trivial. By \Cref{prop:const} and our triviality assumption on $\sE|_{\{t = \infty\}}$, this means that $\sE|_{\bP^1_A \setminus Y}$ is trivial, so that $\sE|_{\{t = 0\}}$ is trivial, too.

To conclude the proof, we now argue that $Y$ as above exists. In fact, it suffices to find an $A$-(finite \'{e}tale) $Y$ as above with the condition $Y \subset \bG_{m,\, A}$ weakened to the condition that there be no finite field obstruction to embedding $Y$ into $\bG_{m,\, A}$: the primitive element theorem for finite separable field extensions will then imply that the embeddings $Y_{k_\fm} \hra \bG_{m,\, k_\fm}$ exist for all maximal ideals $\fm \subset A$ and the Nakayama lemma will   allow us to lift them to an embedding $Y \hra \bG_{m,\, A} \subset \bA^1_A$. To find such a $Y$, we first consider the projective, smooth $A$-scheme $X$ that parametrizes parabolic subgroups of $G$ (see \cite{SGA3IIInew}*{expos\'e XXVI, corollaire~3.5}), so that $X(k_\fm) \neq \emptyset$ for every maximal ideal $\fm \subset A$ with $G_{k_\fm}$ isotropic. The projectivity and smoothness of $X$ allow us to apply the Bertini theorem to iteratively cut $X$ by smooth hypersurfaces passing through specified $k_\fm$-points of $X$ for each  $\fm$ as above to build an $A$-(finite \'{e}tale) $Y_0 = \Spec(A_0) \subset X$ such that $Y_0(k_\fm) \neq \emptyset$ for every maximal ideal $\fm \subset A$ with $G_{k_\fm}$ isotropic; alternatively, to build such a $Y_0$ we may apply~\cite{torsors-regular}*{Lemma~6.2.2}, whose proof gives this Bertini argument in detail. For each $N \ge 1$, consider a finite \'{e}tale cover $Y_N \surjects Y_0$ defined by a monic polynomial $f_N(t) \in A_0[t]$ of degree $N$ whose reduction modulo each maximal ideal $\fn \subset A_0$ is a product of $N$ distinct monic linear factors if $k_\fn$ is infinite (resp.,~is irreducible of degree $N$ if $k_\fn$ is finite). The advantage of $Y_N$ is that there is no finite field obstruction to embedding it into $\bG_{m,\, A}$ granted that $N$ is large, in fact, the same even holds for $Y \ce Y_N \sqcup Y_{N + 1}$. By construction, this $Y$ is as required: $G_Y$ is isotropic (even $G_{Y_0}$ is) and, for each maximal ideal $\fm \subset A$ with $G_{k_\fm}$ isotropic, $Y_{k_\fm}$ has two disjoint clopens of degrees $N$ and $N + 1$ over $k_\fm$. 
\epf
 
 \brems
 \remi
 \Cref{thm:section-const} fails beyond semilocal $A$. Indeed, among the rings of integers $\cO_K$ of number fields $K$ for which the class number is not $1$, one finds plenty of examples of \emph{nonprincipal} ideals $I \subset \cO_K$. Since $I$ is generated by two elements, there exists an $s \in \bP^1_{\cO_K}(\cO_K)$ such that $s^*(\sO(1))$ is isomorphic to $I$ and so is nontrivial. 
 
 \remi
 Even though we do not explicate this, the proof of \Cref{thm:section-const} clearly also generalizes and simplifies the aforementioned \cite{torsors-regular}*{Proposition 5.3.6} (so also \cite{Fed22a}*{Theorem 6}).
 \erems

 
\section{Unramified Grothendieck--Serre in the totally isotropic case} \label{sec:tot-iso}

We are ready to settle the unramified case of the Grothendieck--Serre conjecture for reductive groups whose adjoint quotients are totally isotropic in \Cref{thm:tot-iso} below (see \S\ref{pp:conv} for a review of total isotropicity). The final input to this is a study of torsors over $\bA^1_A$ built on the corresponding study of torsors over $\bP^1_A$ carried out in \S\ref{sec:BunG}. For us, a key advantage of $\bA^1_A$ is that we no longer need to restrict to semilocal $A$ thanks to the following general form of Quillen patching due to Gabber (prior versions \cite{Mos08}*{Satz 3.5.1} or \cite{AHW18}*{Theorem 3.2.5} would also suffice for our purposes). 

\blem[\cite{torsors-regular}*{Corollary 5.1.5}] \label{lem:Quillen}
For a locally finitely presented group algebraic space $G$ over a ring $A$, a $G$-torsor \up{for fppf topology} on $\bA^1_A$ descends to $A$ iff it does so Zariski locally on $\Spec A$.
\elem

The following theorem is our key conclusion about torsors over $\bA^1_A$ and is a positive answer to a generalization of \cite{torsors-regular}*{Conjecture 3.5.1} of Horrocks type. In its statement, even when $A$ is local, we cannot drop total isotropicity, see \cite{Fed16}*{Theorem 3 and what follows}.

\bthm \label{thm:Horrocks}
For a reductive group $G$ over a ring $A$ with $G^{\mathrm{ad}}$ totally isotropic, no nontrivial $G$-torsor over $\bA^1_A$ trivializes over the punctured formal neighborhood $A\llp t\i\rrp$ of the section at infinity\uscolon equivalently, every $G$-torsor $\sE$ over $\bP^1_A$ with $\sE|_{\{t = \infty\}}$ trivial restricts to the trivial torsor over $\bA^1_A$. 
\ethm

\bpf
The two formulations are equivalent due to Henselian invariance and patching for $G$-torsors, see \cite{Hitchin-torsors}*{Theorem 2.1.6 and Lemma 2.2.11~(b)}. Moreover, by base change along the map $\bA^1_A \cong \Spec(A[u]) \ra \Spec A$, we obtain a $G$-torsor $\sE_u$ over $\bP^1_{A[u]}$ with $\sE_u|_{\{t = \infty \}}$ trivial such that the restriction of $\sE_u$ to the ``diagonal'' section $t = u$ of $\bA^1_{A[u]}$ is $\sE$. Thus, by changing the coordinates of $\bP^1_{A[u]}$ via $[x : y] \mapsto [x - uy : y]$ and replacing $A$ and $\sE$ by $A[u]$ and $\sE_u$, respectively, we are left with showing that our $G$-torsor $\sE$ over $\bP^1_A$ with $\sE|_{\{t = \infty\}}$ trivial is such that $\sE|_{\{t = 0\}}$ is also trivial. 

This last claim is insensitive to replacing $\sE$ by its pullback along the map $\varphi_d\colon \bP^1_A \ra \bP^1_A$ given by $[x : y] \mapsto [x^d : y^d]$ for a $d > 0$. We replace $\sE$ by such a pullback with $d$ being the degree of the isogeny $(G^{\der})^{\mathrm{sc}} \times \rad(G) \ra G$. Since the resulting pullback of $\sO(1)$ is $\sO(d)$, by \cite{Gil02}*{th\'{e}or\`eme~3.8}, this ensures that each $\sE|_{\bP^1_{\ov{k}_s}}$ for $s \in S$ now lifts to a $((G^{\der})^{\mathrm{sc}} \times \rad(G))$-torsor over $\bP^1_{\ov{k}_s}$. 

The obtained fibral liftability and \Cref{lem:Cartesian} imply that $\sE$ itself lifts to a $((G^{\der})^{\mathrm{sc}} \times \rad(G))$-torsor over $\bP^1_A$ whose restriction to the section at infinity is trivial, to the effect that we have reduced to $G$ being either a torus or semisimple, simply connected. Moreover, in the toral case, $\sE|_{\bA^1_A}$ is trivial by \Cref{lem:BunM}, so for the rest of proof we assume that $G$ is semisimple, simply connected. Due to \cite{SGA3IIInew}*{expos\'{e} XXIV, section 5.3, proposition 5.10 (i), proposition 8.4} (compare with \eqref{eqn:tot-isot} above), we may then even also assume that $G$ is simple. Granted these reductions, we revert to arguing the triviality of $\sE|_{\bA^1_A}$. For this, we first use \Cref{lem:Quillen} coupled with a limit argument to reduce to the case when $A$ is local. For local $A$, however, $\sE|_{\bA^1_A}$ is trivial by \Cref{cl:isotropic}.
\epf

 We turn to the promised totally isotropic, unramified case of the Grothendieck--Serre conjecture. 
 
 \bthm \label{thm:tot-iso}
 Let $R$ be a Noetherian semilocal ring that is flat and geometrically regular over some Dedekind ring, let $K \ce \Frac(R)$ be its ring of fractions. The Grothendieck--Serre conjecture holds for every  reductive $R$-group $G$ whose adjoint quotient $G^{\mathrm{ad}}$ is totally isotropic, more precisely, for every such $G$, we have
 \[
 \Ker(H^1(R, G) \ra H^1(K, G)) = \{*\}. 
 \]
 \ethm

\bpf
We let $\cO$ be a Dedekind ring over which $R$ is flat and geometrically regular, assume without losing generality that $\cO$ is semilocal, and decompose $\cO$ and $R$ into factors to make them domains. We then combine Popescu's \cite{SP}*{Theorem~\href{https://stacks.math.columbia.edu/tag/07GC}{07GC}} with a limit argument to reduce to when $R$ is the semilocal ring of a smooth, affine, integral $\cO$-scheme $X$. We spread out to ensure that our reductive group $G$ with $G^\ad$ totally isotropic and the generically trivial torsor $E$ that we wish to trivialize both begin life over $X$. 

By \Cref{prop:build-Z} and spreading out, we may replace $X$ by an affine open containing $\Spec R$ to arrange that there be a closed $Z \subset X$ of codimension $\ge 2$ (without loss of generality, cut out by a regular sequence of length $2$---this simplifies the spreading out), an open $V \subset \bP^1_X$ containing both $\bP^1_{X\setminus Z}$ and the $X$-points $\{t = 0 \}$ and $\{t = \infty\}$, and a $G$-torsor $\wt{E}$ over $V$ such that $\wt{E}|_{\{t = 0\}} \simeq E$ and $\wt{E}|_{\{t = \infty\}}$ is trivial. Since $X$ is affine, there is a principal Cartier divisor $Y \subset X$ containing $Z$ and not containing any generic point of any $\cO$-fiber of $X$. Since $X \setminus Y$ is affine, \Cref{thm:Horrocks} ensures that $\wt{E}|_{\bA^1_{X \setminus Y}}$ is trivial, so, by \Cref{thm:Horrocks} again, so is $\wt{E}|_{\bP^1_{X\setminus Y} \setminus \{t = 1\}}$. By patching, then there is a $G$-torsor $\wt{E}'$ over $\bP^1_{X \setminus Y} \cup (V \setminus \{t = 1\})$ that is trivial on $\bP^1_{X \setminus Y}$ and agrees with $\wt{E}$ on $V \setminus \{ t = 1\}$. As in the proof of \Cref{prop:build-Z}, using \cite{CTS79}*{th\'{e}or\`{e}me~6.13} and spreading out, this $\wt{E}'$ extends to a $G$-torsor over $\bP^1_{X \setminus Z'} \cup (V \setminus \{t = 1\}) $ for some closed $Z' \subset Y$ of codimension $\ge 2$ in $X$ containing $Z$. We replace $\wt{E}$ by this extension of $\wt{E}'$ and $Z$ by $Z'$ to assume that our $\wt{E}$ as above trivializes over $\bP^1_{X\setminus Y}$.

If $X$ is of dimension $\le 1$, then $E$ is trivial by \cite{Guo22a}*{Theorem 1}, so we assume that $X$ is of relative dimension $d > 0$ over $\cO$. By \Cref{lem:pres-lem}, we may replace $X$ by an affine open containing $\Spec R$ to find an affine open $S \subset \bA^{d - 1}_\cO$ and a smooth map $X \ra S$ of pure relative dimension $1$ such that $Y \cap X$ is $S$-quasi-finite and $Z \cap X$ is $S$-finite. The base change along $\Spec R \ra S$ then gives
\bitem
\m
a smooth, affine $R$-scheme $C$ of pure relative dimension $1$ equipped with an $s \in C(R)$;

\m
a reductive $C$-group scheme $\sG$ with $s^*(\sG) \cong G$ and a $\sG$-torsor $\sE$ over $C$ with $s^*(\sE) \cong E$;

\m
an $R$-quasi-finite closed $\sY \subset C$ and an $R$-finite closed $\sZ \subset \sY$; and 

\m
a $\sG$-torsor $\wt{\sE}$ over $\bP^1_{C \setminus \sZ}$ such that $\wt{\sE}|_{\{t = 0\}} \simeq \sE|_{C \setminus \sZ}$ and both $\wt{\sE}|_{\{t = \infty\}}$ and $\wt{\sE}|_{\bP^1_{C \setminus \sY}}$ are trivial.
\eitem
As in the proof of \Cref{prop:build-Z}, we will gradually simplify these data to show that $E$ is trivial. The $R$-finiteness of $\sZ$, as opposed to $R$-quasi-finiteness as there, makes some of these simplifications easier, but dragging $\wt{\sE}$ along complicates some others. To begin with, as there, we use \cite{Li24}*{Proposition~7.4} to replace $C$ by a finite \'etale cover of some affine open neighborhood of $\sZ \cup s$ to reduce to when $\sG \cong G_{C}$, compatibly with the identification after $s$-pullback. Similarly, by \cite{split-unramified}*{Lemma~6.1}, we may replace $C$ by a finite \'{e}tale cover of some affine open neighborhood of $\sZ \cup s$ to reduce further to when there is no finite field obstruction to embedding $\sZ \cup s$ into $\bA^1_R$. We then shrink $C$ around $\sZ \cup s$ to ensure that there is no finite field obstruction to embedding $\sY \cup s$ into $\bA^1_R$ either. 

\Cref{lem:embed-field,lem:reembed} now ensure that at the cost of replacing $C$ by an affine open containing the closed $R$-fibers of $\sY \cup s$ (so also containing $\sZ \cup s$), there are an affine open $W \subset \bA^1_R$ and an \'{e}tale $R$-morphism $f\colon C \ra W$ that embeds $\sY \cup s$ excisively into $W$, so that we have a Cartesian square
\[
\xymatrix{
\sY \ar@{=}[d] \ar@{^(->}[r] & C \ar[d]^-f \\
\sY \ar@{^(->}[r] & W
}
\]
in which the horizontal maps are closed immersions. We wish to replace $C$ by $W$, and for this we will now use excision (see \cite{torsors-regular}*{Proposition~4.2.1}) to descend $\wt{\sE}$ to $\bP^1_{W\setminus \sZ}$. First of all, by \Cref{prop:const}~\ref{m:const-a} (by the full faithfulness conclusion applied to the automorphisms of the trivial $G$-torsor), we have $G(C \setminus \sY) \isomto G(\bP^1_{C\setminus \sY})$, so the set of trivializations of $\wt{\sE}|_{\bP^1_{C \setminus \sY}}$ maps bijectively onto its counterpart for $(\wt{\sE}|_{\{t = \infty\}})|_{C\setminus \sY}$. Thus, $\wt{\sE}|_{\bP^1_{C \setminus \sY}}$ has a trivialization $\gA$ whose restriction to the infinity section extends to a trivialization of $\wt{\sE}|_{\{t = \infty\}}$ over all of $C \setminus \sZ$. We use this $\gA$ to descend $\wt{\sE}|_{\bP^1_{C \setminus \sY}}$ to a trivial $G$-torsor over $\bP^1_{W \setminus \sY}$. By excision, the latter then extends uniquely to a $G$-torsor $\wt{\sE}'$ over $\bP^1_{W\setminus \sZ}$ descending $\wt{\sE}$. By excision and the choice of $\gA$, our trivialization of $\wt{\sE}'|_{\bP^1_{W \setminus \sY}}$ restricts to a trivialization of $(\wt{\sE}'|_{\{t = \infty \}})|_{W \setminus \sY}$ that extends to a trivialization of $\wt{\sE}'|_{\{t = \infty \}}$ over all of $W \setminus \sZ$.


At this point we have constructed a $G$-torsor $\sE' := \wt{\sE}'|_{\{t = 0\}}$ over $W \setminus \sZ$ whose base change to $ C \setminus \sZ$ is $\wt{\sE}|_{\{t = 0\}} \simeq\sE|_{C \setminus \sZ}$. However, our \'{e}tale map $f\colon  C \ra W$ is excisive with respect to $\sZ$ as well, so, by excision again, $\sE'$ extends to  a $G$-torsor over all of $W$ that descends $\sE$. We may therefore replace $C$ by $W$ and $\sE$ (resp.,~$\wt{\sE}$) by this extension (resp.,~by $\wt{\sE}'$) to reduce to $C$ being an affine open~of~$\bA^1_R$. 

Once $C$ is an open of $\bA^1_R$, however, the existence of an $R$-point $s$ of $C$ forces $\bP^1_R \setminus C$ to be $R$-finite. The avoidance lemma \cite{GLL15}*{Theorem~5.1} (recalled in \cite{split-unramified}*{Lemma 3.1}) then supplies an $R$-finite hypersurface $H \subset C \subset \bP^1_R$ containing $\sZ$. The complement $C \setminus H$ is affine, so the triviality of $\wt{\sE}|_{\{t = \infty\}}$ and \Cref{thm:Horrocks} ensure that $\wt{\sE}|_{\bA^1_{C\setminus H}}$ and thus also $\sE|_{C \setminus H}$ are trivial. In particular,  since $H$ is closed in $\bP^1_R$, by patching, $\sE$ extends to a $G$-torsor over $\bP^1_R$ that is trivial  at infinity. \Cref{thm:section-const} then ensures that the pullback under $s$, that is, $E$, is trivial as well, as desired. 
\epf

\brem
The proof of \Cref{thm:tot-iso} uses the $G$-torsor $\wt{E}$ over $\bP^1_{X \setminus Z}$ as a ``witness'' of $E$ being simpler  over $X \setminus Z$. At the cost of first passing to simply connected groups via \Cref{prop:pass-to-sc}, one can also carry out the proof with a ``unipotent chain of torsors'' as a witness. Namely, at the cost of shrinking $X$ around $\Spec R$, one may fix sufficiently general opposite proper parabolic subgroups $P^+, P^- \subset G$ and use the Borel--Tits theorem \cite{Gil09}*{fait~4.3, lemme 4.5} (which needs both the total isotropicity and the simply connectedness assumptions) to build a principal Cartier divisor $Y \subset X$, a closed $Z \subset Y$ of codimension $\ge 2$ in $X$, and a sequence $E_0, \dotsc, E_n$ of $G$-torsors over $X \setminus Z$ such that 
\bitem
\m
each $E_i$ is trivialized over $X \setminus Y$, the $(X \setminus Z)$-group $\Aut_G(E_i)$ has opposite parabolic subgroups $P^{\pm}_i$ that under the trivialization over $X \setminus Y$ correspond to $P^{\pm}|_{X\setminus Y}$, and the $\Aut_G(E_i)$-torsor $\Isom_G(E_i, E_{i + 1})$ for $i < n$ reduces either to a $\sR_u(P_i^+)$-torsor or to a $\sR_u(P_i^-)$-torsor over~$X \setminus Z$;

\m
$E_0$ is trivial and $E_n$ is the restriction of our generically trivial $G$-torsor $E$ over $X$ to $X \setminus Z$.
\eitem
Since torsors under unipotent radicals of parabolic subgroups trivialize over affine schemes (see \cite{SGA3IIInew}*{expos\'{e} XXVI, corollaire 2.5}), the existence of the ``unipotent chain'' $E_0, \dotsc, E_n$ implies that $E$ trivializes over every affine $(X \setminus Z)$-scheme, and it is possible to carry out the proof of \Cref{thm:tot-iso} by dragging the chain $E_0, \dotsc, E_n$ along in place of $\wt{E}$ in the intermediate steps. 

For a systematic development of the notion of a unipotent chain of torsors, see \cite{Fed23}.
\erem

 
 \section{Reducing to semisimple, simply connected groups} \label{sec:reduce-ss}

 We combine the work of \S\S\ref{sec:codim-2}--\ref{sec:BunG} with purity theorems for $H^{\le 2}$ with multiplicative group coefficients (essentially, purity for the Brauer group \cite{brauer-purity}) to reduce the unramified case of the Grothendieck--Serre conjecture to simply connected $G$. The method is new even in equal characteristic, although the corresponding reduction in equal characteristic was  the main goal of the article \cite{Pan20b}.

 \bprop \label{prop:pass-to-sc}
Let $G$ be a reductive group over a Noetherian semilocal ring $R$ that is flat and geometrically regular over some Dedekind ring. Every generically trivial $G$-torsor over $R$ lifts to a generically trivial $(G^\der)^{\mathrm{sc}}$-torsor over $R$ \up{with notation as in \uS\uref{pp:conv}}, so, setting $K \ce \Frac(R)$,~we~have
\[
\Ker(H^1(R, (G^\der)^{\mathrm{sc}}) \ra H^1(K, (G^\der)^{\mathrm{sc}})) = \{ *\} \q \Longrightarrow \q \Ker(H^1(R, G) \ra H^1(K, G)) = \{ *\}.
\]
\eprop

\bpf
For a generically trivial $G$-torsor $E$ over $R$ to be lifted to a generically trivial $(G^\der)^{\mathrm{sc}}$-torsor, \Cref{prop:build-Z} gives us an open $V \subset \bP^1_R$ containing $\{t = 0\}$ and $\{t = \infty\}$ with complement $\bP^1_R \setminus V$ of codimension $\ge 3$ in $\bP^1_R$ and a $G$-torsor $\sE$ over $V$ such that $\sE|_{\{t = 0\}} \simeq E$ and $\sE|_{\{t = \infty\}}$ is trivial. It suffices to lift some twist of $\sE$ by an $R$-sectionwise trivial $\Rad(G)$-torsor over $V$ to a $(G^\der)^{\mathrm{sc}}$-torsor $\wt{\sE}$ over $V$ with $\wt{\sE}|_{\{t = \infty\}}$ trivial: then $\wt{\sE}|_{\{t = 0\}}$ will lift $E$ and be generically trivial by \Cref{thm:section-const} applied with $A = K$.

Set $Z \ce \Ker((G^\der)^{\mathrm{sc}} \ra G)$. By the codimension condition and purity \cite{flat-purity}*{Theorem 7.2.9},
\be \label{eqn:purity}
H^1(\bP^1_R, \Corad(G)) \isomto H^1(V, \Corad(G)) \qxq{and} H^2(\bP^1_R, Z) \isomto H^2(V, Z).
\ee
In particular, the $\Corad(G)$-torsor induced by $\sE$ extends to a $\Corad(G)$-torsor over $\bP^1_R$ that is trivial at infinity and hence, by \Cref{lem:BunM}, comes from $\sO(1)$ via a cocharacter $\bG_{m,\, R} \ra \Corad(G)$. Thus, since $\Rad(G) \ra \Corad(G)$ is an isogeny, as in the proof of \Cref{thm:section-const}, by pulling back along the base change to $V$ of the map $\varphi_d \colon \bP^1_R \ra \bP^1_R$ for some $d > 0$ such that $\varphi_d$ sends the homogeneous coordinates of $\bP^1_R$ to their $d$-th powers, we reduce to the case when the $\Corad(G)$-torsor induced by $\sE$ lifts to an $R$-sectionwise trivial $\Rad(G)$-torsor. By twisting $\sE$ by such a lift, we may assume that $\sE$ induces a trivial $\Corad(G)$-torsor, so lifts to $G^\der$-torsor over $V$. By \cite{Gir71}*{chapitre~III, proposition~3.3.3~(iv)}, the group $\Corad(G)(V)$  acts transitively on the set of isomorphism classes of such lifts over $V$, and likewise after restricting to the infinity section. Thus, since this restriction induces a surjection $\Corad(G)(V) \surjects \Corad(G)(R)$, we may lift $\sE$ to a $G^\der$-torsor whose restriction to infinity is trivial. In effect, we may replace $G$ by $G^\der$ to reduce to the case when $G$ is semisimple.  

Once $G$ is semisimple, the obstruction to lifting $\sE$ to a $G^{\mathrm{sc}}$-torsor lies in $H^2(V, Z) \cong H^2(\bP^1_R, Z)$. By replacing $V$ by its pullback by $\varphi_d$ for some $d > 0$ and applying \cite{Gil02}*{th\'eor\`eme 3.8} as in the proof of \Cref{thm:section-const}, we may arrange that the restriction $\sE|_{\bP^1_{\ov{K}}}$ to the geometric generic fiber lifts to a $G^{\mathrm{sc}}$-torsor over $\bP^1_{\ov{K}}$, in other words, that the obstruction in question vanishes after pullback to $\bP^1_{\ov{K}}$. By the triviality at infinity and \Cref{lem:gerbeM}, however, it then vanishes already over $V$, to the effect that $\sE$ lifts to a $G^{\mathrm{sc}}$-torsor over $V$. By \cite{Gir71}*{chapitre~III, proposition~3.4.5~(iv)}, the group $H^1(V, Z)$ acts transitively on the set of isomorphism classes of such lifts. Thus, since restriction to infinity induces a surjection $H^1(V, Z) \surjects H^1(R, Z)$, a desired lift $\wt{\sE}$ indeed exists.
\epf






\end{document}